\documentclass[12pt,letterpaper]{amsart}
\usepackage{amsthm,amssymb,latexsym,amsmath}
\usepackage{graphicx,color,xcolor,comment}

\usepackage{booktabs}
\usepackage{graphicx}
\usepackage{url}
\usepackage{xcolor}
\usepackage{hyperref}
\usepackage{natbib}

\title[Modeling ultrasound  impact on  cancer stem cells]{Modeling low-intensity ultrasound mechanotherapy impact on growing cancer stem cells}

\author[B.~Blanco {\it et al}]{B. Blanco}%
\address{B.~Blanco, Department of Electronics and Computing, University of Santiago de Compostela, Spain \& Instituto de Investigación Biosanitaria, Research Unit ``Modelling Nature'' (MNat), University of Gra\-nada, Spain}
\author[]{R. Palma}
\address{R. Palma, Department of Structural Mechanics, University of Gra\-nada, Spain}
\author[]{M. Hurtado}
\address{M. Hurtado, G. Rus, Instituto de Investigaci\'on Biosanitaria, Research Unit ``Modelling Nature'' (MNat) \& Department of Structural Mechanics, University of Granada, Spain}
\author[]{G. Jim\'enez, C. Gri\~n\'an-Lis\'on}
\address{G. Jim\'enez, C. Gri\~n\'an-Lis\'on, J. A. Marchal, Instituto de Investigaci\'on Biosanitaria, Research Unit ``Modelling Nature'' (MNat) \& Department of Human Anatomy and Embryology, University of Granada, Spain}
\author[]{J. Melchor}
\address{J. Melchor, Instituto de Investigaci\'on Biosanitaria, Research Unit ``Modelling Nature'' (MNat) \& Department of Statistics and Operations Research, University of Granada, Spain}
\author[]{J.A. Marchal}
\author[]{H. Gomez}
\address{H. Gomez, School of Mechanical Engineering, Weldon School of Biomedical Engineering, \& Purdue Center for Cancer Research, Purdue University, West Lafayette, Indiana, United States of America}
\author[]{G. Rus}
\author[]{J. Soler}
\address{J. Soler, Department of Applied Mathematics \& Research Unit ``Modelling Nature'' (MNat), University of Granada, Spain}

\begin{document}
\maketitle

\begin{abstract}
Targeted therapeutic interventions utilizing low-inten\-sity ultrasound (LIUS) exhibit substantial potential for hindering the proliferation of cancer stem cells. This investigation introduces a multiscale model and computational framework to comprehensively explore the therapeutic LIUS on poroelastic tumor dynamics, thereby unraveling the intricacies of mechanotransduction mechanisms at play. Our model includes both macroscopic timescales encompassing days and rapid timescales spanning from microseconds to seconds, facilitating an in-depth comprehension of tumor behavior. We unveil the discerning suppression or reorientation of cancer cell proliferation and migration, enhancing a notable redistribution of cellular phases and stresses within the tumor microenvironment.
Our findings defy existing paradigms by elucidating the impact of LIUS on cancer stem cell behavior. This endeavor advances our fundamental understanding of mechanotransduction phenomena in the context of LIUS therapy, thus underscoring its promising as a targeted therapeutic modality for cancer treatment.
Furthermore, our results make a substantial contribution to the broader scientific community by shedding light on the intricate interplay between mechanical forces, cellular responses, and the spatiotemporal evolution of tumors. These insights hold the promising to promote a new perspective for the future development of pioneering and highly efficacious therapeutic strategies for combating cancer in a personalized manner.
\end{abstract}



\section{Introduction}\label{section1}
Mechanotherapy represents an emerging frontier in cancer treatment, harnessing the power of mechanical forces {\it per se}, or through their interaction with cell biochemical connections, to selectively target and eradicate or reverse the growth trend of cancer cells. The underlying principle behind this strategy lies in the observation that cancer cells exhibit increased sensitivity to mechanical stimuli, and by manipulating these forces, their properties can be affected, ultimately leading to cell death or dysfunction.

In the past few years, extensive endeavors have been dedicated to develop a wide range of techniques aimed at modulating cell behavior by altering the microenvironment. These techniques encompass a spectrum of perspectives, from pharmacological agents that modify the elasticity of the remodeled microenvironment and cell stiffness~\citep{huang2020stiffness,jain2010delivering,polydorou2017pirfenidone,panagi2022polymeric,abedi2022ultrasound}, to the application of mechanical waves. These principles not only serve as fundamental research tools for investigating the basic interactions of cell mechanics, but an extraordinary translational potential could also emerge for clinical applications, since mechanical waves can be delivered to patients using transducers or patches. In addition, there is growing recognition that the combination of these therapies can synergistically enhance treatment effectiveness and greatly improve the overall prognosis of the disease.

Recent studies have provided compelling evidence of the therapeutic effectiveness of low-intensity ultrasound. Low-intensity ultrasound (LIUS) and its pulsed version (LIPUS), have been proposed to impact cancer cells by two main mechanisms: i) selectively resonating the right diameter cells under the name of \textit{oncotripsy}, which lies on destroying the cytoskeleton via cavitation~\citep{mittelstein2020selective,heyden2016oncotripsy,heyden2017investigation,lin2022low,lucchetti2020low,prentice2005membrane}, and ii) triggered response produced via mechanotransduction signaling pathways~\citep{lin2022low,lucchetti2020low,katiyar2020inhibition,carina2018inhibitory,gonzalez2023low,singh2021enhanced,tijore2020ultrasound,na2008rapid,geiger2009environmental,vogel2006mechanotransduction,blanco2023review,broders2018mechanotransduction}. 

Although these studies have repeatedly evidenced considerable promising effects, the lack of understanding of the mechanism, and the even opposing responses triggered by diverse frequencies, energies, and configurations, make the concept of little use at its current state. 

The deepening in the configuration and effects of the application of mechanical waves have garnered considerable attention in the scientific community, driven by the potential of LIUS to influence cancer cell dynamics. However, our understanding of the underlying mechanisms behind LIUS-induced effects in cancer cells has been hampered by the complexity of accounting for the combined effects of the large number of agents involved in this process, and the enormous costs associated with performing extensive biological experiments.

To unravel the intricate mechanisms of LIUS and improve its therapeutic prospective, mathematical oncology emerges as a valuable tool \citep{agus2012sciences}. By simulating the complex interactions between ultrasound waves and the tumor microenvironment, in conjunction with other treatment modalities, these models offer a great avenue to try to understand LIUS at a deeper level. Such ultrasound-tumor interactions provide crucial insights into the underlying mechanisms of LIUS and can facilitate the development of more efficient treatment strategies. 

Within this framework, we present a multiscale model to unveil the influence of LIUS on tumor evolution through mechanotransduction. Our approach encompasses the application of LIUS at high frequencies coupled with lower acoustic pressures to target cancer precursor cells known as cancer stem cells (CSCs)~\citep{colak2014cancer,olivares2020csc}, which are widely considered to have an important impact on cancer metastasis and are frequently associated with relapse due to their self-renewal, differentiation capabilities and resistance to conventional therapies~\citep{olivares2020csc}.

The multiscale mathematical model that we present in Section~\ref{section2} characterizes tumors as poroelastic materials composed of an interstitial fluid phase and distinct solid phases governed by elastic properties attributed to tumor cells, healthy cells, and the extracellular matrix (ECM). We coupled the influence of ultrasound to this system and introduced a new mechanotransduction function sensitive to hydrostatic stress. 

\section{Materials and methods}\label{section2}
\subsection{Mathematical model}
The multiscale model is proposed on two different scales: i) slow time scale, $t$, in which the tumor grows and migrates, and ii) fast time scale, $t_u$, in which ultrasound propagates through the tumor. Both scales are coupled by mechanotransduction, which occurs at an ultrasonic time interval and triggers a tumor dynamics response at a slow scale. Thus, we can encapsulate the enduring effects of growth and reorganization, which may not be discernible on the ultrasonic scale but acquire significance on a broader and slower scale~\citep{rus2014nature}.

We describe tumors as poroelastic materials composed of a fluid phase ($\phi_F$) of interstitial fluid and different solid phases ($\phi_i$) that provide elastic stiffness. The solid phases included here are tumor cells ($\phi_T$), healthy cells ($\phi_H$), and the extracellular matrix ($\phi_M$). 
The proposed system is built upon Biot's poroelasticity and growth theory, widely studied in the field of thermodynamics. For a comprehensive grasp of the equations, particularly those pertaining to poroelastic cell competition and mechanotransduction, we suggest consulting references~\citep{biot1941general,carotenuto2018growth,carotenuto2021lyapunov,blanco2023review}. We consider infinitesimal strain theory and linear elasticity assuming that there is no large deformation during ultrasound insonation~\citep{mittelstein2020selective} and growth does not develop great deformation. Furthermore,  we use the $\pmb{u}-p$ poroelastic notation, neglecting the relative fluid-solid displacement.

The momentum balance that describes the dynamic mechanical equilibrium is:
\begin{equation}
\nabla \cdot \pmb{\sigma} =  \rho \frac{\partial ^2 \pmb{\mathrm{u}}}{\partial {t}^{2}},
\label{balancemultiscale}
\end{equation}
where $\rho$ is the medium density and $\pmb{\mathrm{u}}$ are the displacements. The multiscale Cauchy stress tensor in a sonicated growing tumor is described by:
\begin{equation}
\pmb{\sigma}(\pmb{x},t) = \pmb{\sigma}_s (\pmb{x},t) + \pmb{\sigma}_{u}(\pmb{x},t_{u}),
\label{multiscalestress}
\end{equation}
where the slow-scale stress $\pmb{\sigma}_s (\pmb{x},t)$ accounts for the growth and the poroelastic rearrangements while the fast-scale stress $\pmb{\sigma}_u(\pmb{x},t_{u})$ is the ultrasonic stress. 
To isolate the governing equations at each temporal scale and save computational cost, we follow the principles of multiscale developed in~\citep{kevorkian2012multiple} and adapted to the ultrasound formulation in~\citep{rus2014nature}. Then, we define the average of the multiscale stress over an ultrasonic spatial and temporal cycle, specifically, the reference ultrasonic wavelength $\lambda$ and period $T$:
\begin{equation}
\langle \pmb{\sigma} \rangle =\frac{1}{\lambda T} \int_{0 }^{T} \int_{0 }^{ \lambda} \pmb{\sigma} \,d\pmb{x}\,dt_{u}.
\end{equation}
Considering the definition of multiscale stress, the above reads:
\begin{equation}
\langle \pmb{\sigma} \rangle = \frac{1}{\lambda T} \int_{0 }^{T} \int_{0 }^{ \lambda} \pmb{\sigma}_s \,d\pmb{x}\,dt_{u} + \frac{1}{\lambda T} \int_{0 }^{T} \int_{0 }^{ \lambda} \pmb{\sigma}_u \,d\pmb{x}\,dt_{u},
\label{promediastress}
\end{equation}
where the slow stress independent of the ultrasonic scale is $ \langle \pmb{\sigma}_s\rangle=\pmb{\sigma}_s$. Ultrasonic stress is a sinus function on $\lambda$ and $T$, so $\langle \pmb{\sigma}_u\rangle=0$. Then, the average of the multiscale stress is the slow-scale stress, $\langle \pmb{\sigma} \rangle=\pmb{\sigma}_s$, and subsequently the average of the slow-scale stress is the total multiscale stress, $\langle \pmb{\sigma}_s\rangle=\pmb{\sigma}$. Finally, with respect to equation~\eqref{multiscalestress} and the independence of the slow-scale stress from the ultrasonic scale, the ultrasonic stress is $\pmb{\sigma}_u = \pmb{\sigma} - \langle \pmb{\sigma} \rangle$. 

\begin{figure*}[h]
\centering
\includegraphics[width=\textwidth]{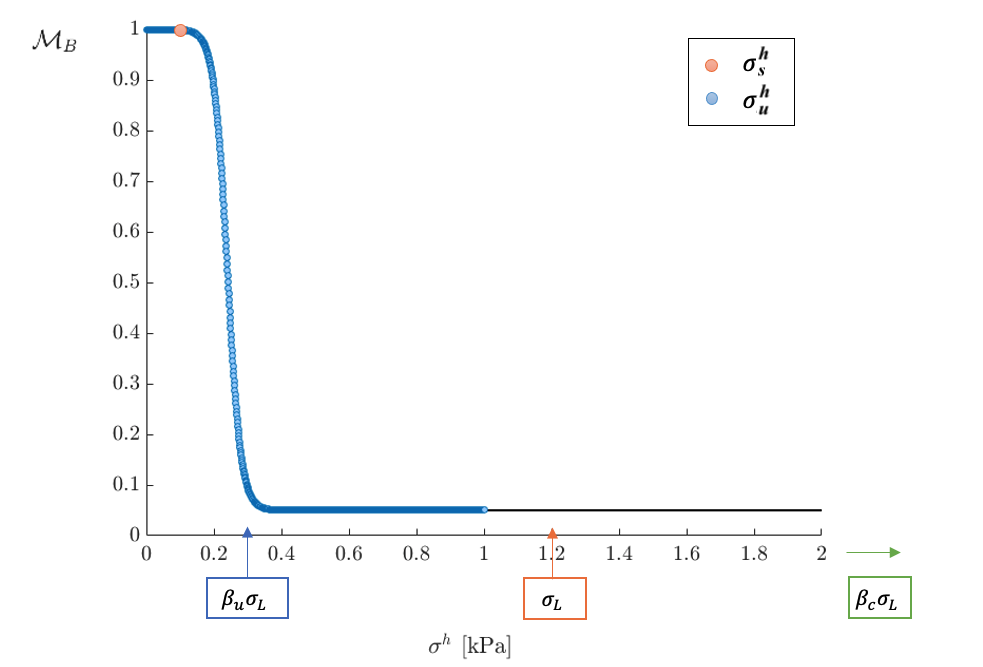}
\caption[Mechanotransduction function]{\textbf{Mechanotransduction function in an ultrasonic interval}. Cells perceive the average of the sigmoid function $\mathcal{M}_{B}$. The slow ultrasound stress is constant at an ultrasonic time interval, while ultrasound stress exhibits dynamic behavior, oscillating between rarefaction and compression -- for this case, we have plotted the absolute stress of a wave with an amplitude of 1kPa --. To account for the dynamic nature of ultrasound stress, the static stress limit $\sigma_L$ is decreased by a coefficient $\beta_u$. However, if the limit is exceeded, it may result in cell disruption and cessation of proliferation or migration, indicated by $\mathcal{M}_B=0$.}
\label{mechanoexp}
\end{figure*}

 Once the multiscale approach is formalized, we define the slow-scale stress as an additive decomposition:
\begin{equation}
\pmb{\sigma}_s = \pmb{\sigma}_{e} + \pmb{\sigma}_{p} + \pmb{\sigma}_{g},
\label{total}
\end{equation}
where $\mathrm{\pmb{\sigma}_{e}}$ is the so-called effective solid stress tensor, $\mathrm{\pmb{\sigma}_{p}}$ the fluid pressure contribution, and $ \mathrm{\pmb{\sigma}_{g}} $  the stress generated during growth. Hence, the equation of equilibrium~\eqref{balancemultiscale} applied to slow-scale stress can be considered as a quasistatic process since characteristic velocities are small and inertia terms can be neglected~\citep{lorenzo2019}. For an elastic isotropic material, the constitutive equation for the effective solid stress that accounts for the elastic rearrangements yields:
\begin{equation}
\pmb{\sigma}_{e} = 2\mu_d\Big(\pmb{\varepsilon}-\dfrac{1}{3}\mathrm{tr}(\pmb{\varepsilon})\pmb{\mathrm{I}} \Big) + K_d  \mathrm{tr}(\pmb{\varepsilon})\pmb{\mathrm{I}},
\label{elas}
\end{equation}
where the small strain is $\pmb{\mathrm{\varepsilon}}=\dfrac{1}{2}(\nabla \boldsymbol{\mathrm{u}} + \nabla \boldsymbol{\mathrm{u}}^\mathrm{T})$, with $\pmb{\mathrm{u}}$ the displacements, $\pmb{\mathrm{I}}$ the second-order identity tensor, and $K_d$ and $\mu_d$ the drained bulk and shear modulus. We can neglect the viscous solid contribution in the slow-scale governing equation since the relaxation terms of rearrangements are on a smaller time scale than growth. 
\noindent The stress produced by the fluid is:
\begin{equation}
\pmb{\sigma}_{p} =-\alpha \Big(p -p_0 \Big)\pmb{\mathrm{I}},
\label{press}
\end{equation}
with $\alpha$ the Biot coefficient, $p$ the fluid pore pressure, and $p_0$ the initial fluid pore pressure. The evolution of the fluid pressure $p$ is regulated by the storage equation:
\begin{equation}
\frac{\partial \zeta}{\partial t} =  \frac{1}{M} \frac{\partial p}{\partial t} + \alpha \frac{\partial tr(\pmb{\varepsilon})}{\partial t}  = \nabla \cdot \Big(k\nabla p\Big) + \mathrm{\Gamma_F},
 \label{porogrowth}
\end{equation}
where $\zeta$ is the dimensionless variation of fluid content defined by the difference between the actual and initial fluid phase $\mathrm{\zeta=\phi_F-\phi_{F0}}$.The parameter $M$ represents the Biot modulus, and $k$ denotes the hydraulic conductivity, given by $k= \kappa {\nu_f }^{-1}$, where $\kappa$ stands for the permeability of the medium, and the dynamic fluid viscosity is described by $\nu_f$.  The source term $\Gamma_F$ accounts for the fluid interchange between vessels and capillaries. Considering the theory of Starling~\citep{carotenuto2021lyapunov,carotenuto2018growth,wu2013effect,fraldi2018cells,stylianopoulos2013coevolution}, the fluid flow source yields:
\begin{equation}
 \Gamma_F = k_v \Big [ (p_v - p) - \omega( \pi_v -\pi_l) \Big] - k_l (p -p_l),
 \end{equation}
where $p_v$ is associated with the vessel pressure, $\omega$ represents the reflection coefficient, weighing the interstitial osmotic pressure $(\pi_v - \pi_l)$, and $p_l$ denotes the lymphatic pressure drainage operating counter to the vessel pressure system. The constants $k_v$ and $k_l$ correspond to the conductivity coefficients of the vessel and lymphatic system, respectively. Following recent literature, we formulate the conductivity of the lymphatic system as a function of tumor cells, encompassing the diminishing drainage of the lymphatic system induced by tumor growth.
\begin{equation}
k_l = \Big [ 1 - (\phi_T - \phi_{T0})\Big]k_{ln},
\end{equation}
where $k_{ln}$ is the conductivity of the lymphatic system under normal conditions~\citep{carotenuto2021lyapunov,wu2013effect}. Finally, the stress produced by growth reads:
\begin{equation}
\pmb{\sigma}_{g} =-K g \pmb{\gamma} ,
\label{sigma_press}
\end{equation}
where $g$ is the growth strain function and $\pmb{\gamma}$ is the tensor that distributes the growth in different directions. In this study, we have assumed isotropic growth, so $\pmb{\gamma}= \frac 1 3 \pmb{\mathrm{I}}$. The growth function is considered homogeneous and, therefore, can be written as:
\begin{equation}
g = \phi_{T}+\phi_{H}+\phi_{M} - \phi_{T0}-\phi_{H0}-\phi_{M0},
 \label{growth}
\end{equation}
with the zero subindexes denoting the initial volume fractions. The volumetric fractions evolve and interact with the mechanical environment and are governed by: 
\begin{equation}
\begin{split}
\frac{\partial \phi_{T}}{\partial t} &= \nabla \cdot \left( \mathcal{M}_T D_T \phi_{T}\nabla \phi_{T}\right)+ \mathrm{\phi_{F} \mathcal{M}_T  \phi_{T} \Gamma_T T_T},\\
\frac{\partial \phi_{H}}{\partial t} &= \phi_{F} \mathcal{M}_H  \phi_{H} \Gamma_H T_H,\\
\frac{\partial \phi_{M}}{\partial t} &= \beta_{T} \phi_{T} + \beta_{H} \phi_{H} -\delta_{M} \phi_{M}\Gamma_{M},
\end{split}
\end{equation}
where the first equation describes the tumor cell dynamics. In particular, the first term on the right hand side accounts for tumor non-linear cell flux, described here by a finite speed tumor propagation front limited by the diffusion coefficient $D_T$, although 
 a controlled velocity of propagation could be also taken into account~\citep{blanco2023review,conte2021modeling,blanco2021modeling,calvo2017qualitative,calvo2016pattern}. The second term considers the competitive interaction among other species -- $\mathrm{\Gamma_T}$ --, and both terms account for the mechanotransduction function -- $\mathrm{ \mathcal{M}_T}$ --. Mechanotransduction and competition are also described for healthy cells by $\mathrm{\mathcal{M}_H}$ and $\mathrm{\Gamma_H}$. The ECM evolution depends on the species interaction $\mathrm{\Gamma_M}$ and on the ECM synthesis promoted by the cells by the production rates $\beta_T$ and $\beta_H$, and the ECM degradation processes enabled by the loss rate $\delta_M$~\citep{carotenuto2021lyapunov}. The competition terms $\Gamma_T$, $\Gamma_H$ and  $\Gamma_M$ are defined by the following Volterra-Lokta dynamics, see~\citep{carotenuto2021lyapunov}:
\begin{equation}
\begin{split}
 \Gamma_T &= \Big( 1 - \alpha_{TT}\phi_{T} -\alpha_{TH}\phi_{H} -\alpha_{TM}\phi_{M}\Big),\\
\Gamma_H &= \Big( 1 - \alpha_{HT} \phi_{T} -\alpha_{HH}\phi_{H} -\alpha_{HM}\phi_{M}\Big),\\
\Gamma_M &=\alpha_{MT}\phi_{T} + \alpha_{MH}\phi_{H},
\end{split}
\end{equation}
where the coefficients $\alpha_{ij}$, with $i,j=\{T,H,M\}$,  represent the interaction among the cell species.  
To complete the system of equations, we define the mechanotransduction function based on previously validated expressions~\citep{blanco2023review}.  Thus, the function of mechanotransduction is defined in an ultrasonic time period in which cell mechanosensors could receive signaling linked to the cytoskeleton network extremely quickly~\citep{hoffman2011dynamic,na2008rapid}, and then respond triggering changes in proliferation and migration, as proposed in~\citep{geiger2009environmental,hoffman2011dynamic}. Then, we propose that cells perceive the average of the sigmoid function $ \mathcal{M}_{B_i}$ at an ultrasonic time interval, in which cells could sense perturbations and activate mechanotransduction pathways that alter proliferation above a certain stress threshold~\citep{vogel2006mechanotransduction,broders2018mechanotransduction} -- see Figure~\ref{mechanoexp} --. Then, mechanotransduction can be expressed as:
\begin{equation}
\mathcal{M}_i =   \dfrac{1}{T} \int_{0}^{T} { \mathcal{M}_{B_i} dt_u},
\label{mechano2}
\end{equation}
where $\mathcal{M}_{B_i}$ is based on~\citep{carotenuto2021lyapunov}:
\begin{equation}
\begin{split}
\mathcal{M}_{B_i} = \Big[q_i + (1- q_i ) e^{b_i (|\sigma_{\mathcal{M}}| - \beta_s\sigma_{L_i})}\Big] \\ \Big(1 + e^{b_i(|\sigma_{\mathcal{M}}| - \beta_s\sigma_{L_i}) } \Big)^{-1}.
\end{split}
\end{equation}
Indeed, the initial proliferation or migration of cells decreases to the maximum of the viability of the cells, achieving the factor of $q_i$ when the stress perceived by the cells $\sigma_{\mathcal{M}}$ in the environment exceeds a threshold $\sigma_{L_i}$. In the literature, this threshold is obtained for static stress and values between $[1-10]$kPa~\citep{helmlinger1997,roose2003solid,cheng2009,carotenuto2021lyapunov,carotenuto2018growth,fraldi2018cells}. However, to also account for dynamic stress, we adopt a linear parameter $\beta_u$, which reduces the sensitivity limit of cells. This hypothesis is rooted in the understanding that static stress necessitates higher intensity to elicit a response due to stress dissipation, while dynamic pressure operates within a compressed time frame, precluding dissipation. Key factors influencing dissipation include the cytoskeleton, which imparts structural integrity to cells and facilitates the redistribution of mechanical loads within them, as well as the viscosity of solid phases and the dissipation of stresses through interstitial fluid perfusion via pores. The parameter $b_i$ refers to the smoothness of the transition zone of the sigmoid function and determines how fast or slow cells adapt their proliferation to stress. 

Then, cells detect both static hydrostatic growth-induced stress and dynamic ultrasonic-induced stress through mechanotransduction pathways and the total stress perceived yields: 
\begin{equation}
\sigma_{\mathcal{M}}=\sigma^h_s +\sigma^h_{u},
\label{sigmaperc}
\end{equation}
where the superscript $h$ denotes the hydrostatic stress defined by $\sigma^h = \frac{1}{3} \mathrm{tr}( \pmb{\sigma})$ for each time-scale stress. Shear stress contribution is disregarded in our analysis due to the plane ultrasound wave and isotropic growth, as its magnitude is three orders of magnitude smaller than that of the normal components. 

We obtain $\pmb{\sigma}_u$ over a period of time from the propagation of a P-wave emitted by a transducer through the medium. We formulate the displacement of the wave as a boundary condition, for instance, in a lateral face. Then, displacements generated by a transducer in the y-axis direction are described in the simplified form:
\begin{equation}
\pmb{\mathrm{u}}_u=\Big(0 , A \sin{(2\pi{}f t_u)}\Big),
\label{propagation}
\end{equation}
where $A$ is the wave amplitude in terms of displacements and $f$ is the central frequency.  Before natural attenuation, the wave travels at speed $c_p=\sqrt{\frac{K + 4/3\mu}{\rho}}$, where $K$ and $\mu$ are the undrained bulk and shear modulus, and $\rho$ stands for medium density.

The dynamic balance equation accounts for the inertial terms produced during sonication can be written as:
 \begin{equation}
\nabla \cdot \pmb{\sigma}_u= \rho \frac{\partial ^2 \pmb{\mathrm{u}}_u}{\partial {t_{u}}^{2}},
\label{BMD}
\end{equation}
where $\mathrm{\pmb{\sigma}_{u}}$ is the stress produced by mechanical wave propagation. To accurately represent the complex attenuation that occurs on a fast time scale, we utilize the Kelvin-Voigt governing equation as presented below:
 \begin{equation}
\pmb{\sigma}_u= 2\mu \Big(\pmb{\varepsilon}_u-\dfrac{1}{3}\mathrm{tr}(\pmb{\varepsilon}_u)\pmb{\mathrm{I}} \Big) + K \mathrm{tr}(\pmb{\varepsilon}_u)\pmb{\mathrm{I}} + \eta_{K} \frac{\partial \mathrm{tr}(\pmb{\varepsilon}_u)}{\partial t_u} \pmb{\mathrm{I}},
 \label{sigultrasound}
\end{equation}
where the small strain is $\pmb{\varepsilon}_u=\dfrac{1}{2}(\nabla \boldsymbol{\mathrm{u}}_u + \nabla {\boldsymbol{\mathrm{u}}_u}^\mathrm{T})$. Attenuation is described by the volumetric viscosity $\eta_{K}$, neglecting the contribution of the shear viscosity due to the low order of magnitude of the shear component of the compression waves. As described in previous works~\citep{dukhin2009bulk,claes2021measurement},  we define $\eta_K = \frac{\alpha_{\eta} 2 \rho {c_p}^3}{(2\pi f)^2}$, considering the attenuation coefficient $\alpha_{\eta}$ of an ultrasonic wave at a given frequency~\citep{d1986frequency}.

\subsection{Numerical Methods}

Regarding the initial conditions, the initial fluid phase is defined by the equation $\mathrm{\phi_{F0}=1-\phi_{T0}-\phi_{H0}-\phi_{M0}}$. The initial fluid pressure guarantees the equilibrium of the storage equation, so $p_{0}$ causes the source term to be null at the initial time instant $\mathrm{t=0}$h. Furthermore, the initial components of the tumor and healthy cells are distributed in space according to a smoothing function $S$: 
\begin{equation}
\begin{split}
\phi_{T0}&=\overline{\phi}_{T0} S,\\
\phi_{H0}&=\overline{\phi}_{H0} \Big(1-S\Big),\\
S &= \Big[ 1 + e^{b_S \frac{(r-l_t)}{l}}\Big]^{-1},
\end{split}
\end{equation}
where the parameter $\overline{\phi}_{i0}$ is the initial concentration rate,  $l_t$ represents the tumor size, $r$ the radial coordinate, and $l$ the total length of the medium, while $b_S$ is the smoothing coefficient, according to reference~\citep{carotenuto2021lyapunov}.

Regarding boundary conditions (BC), we adopt Winkler-inspired boundary conditions to consider tumor spheroid confinement at slow scales, ~$\pmb{\mathrm{\sigma}_n}=-k_w \pmb{\mathrm{u}}$, where $\pmb{\mathrm{n}}$ is the outer normal vector and $k_w$ is a constant~\citep{lorenzo2019}.  For fast-scale ultrasound propagation, we use the Lysmer-Kuhlemyer boundary condition to account for non-reflecting boundaries~\citep{lysmer1969finite}. The normal stress reads $\pmb{\mathrm{\sigma}_n} =k_a \rho c_p \frac{\partial  \pmb{\mathrm{u}}}{\partial t}$, where $k_a$ is a constant in the range $[0-100]$.

\begin{figure*}[h]
\centering
\includegraphics[width=\textwidth]{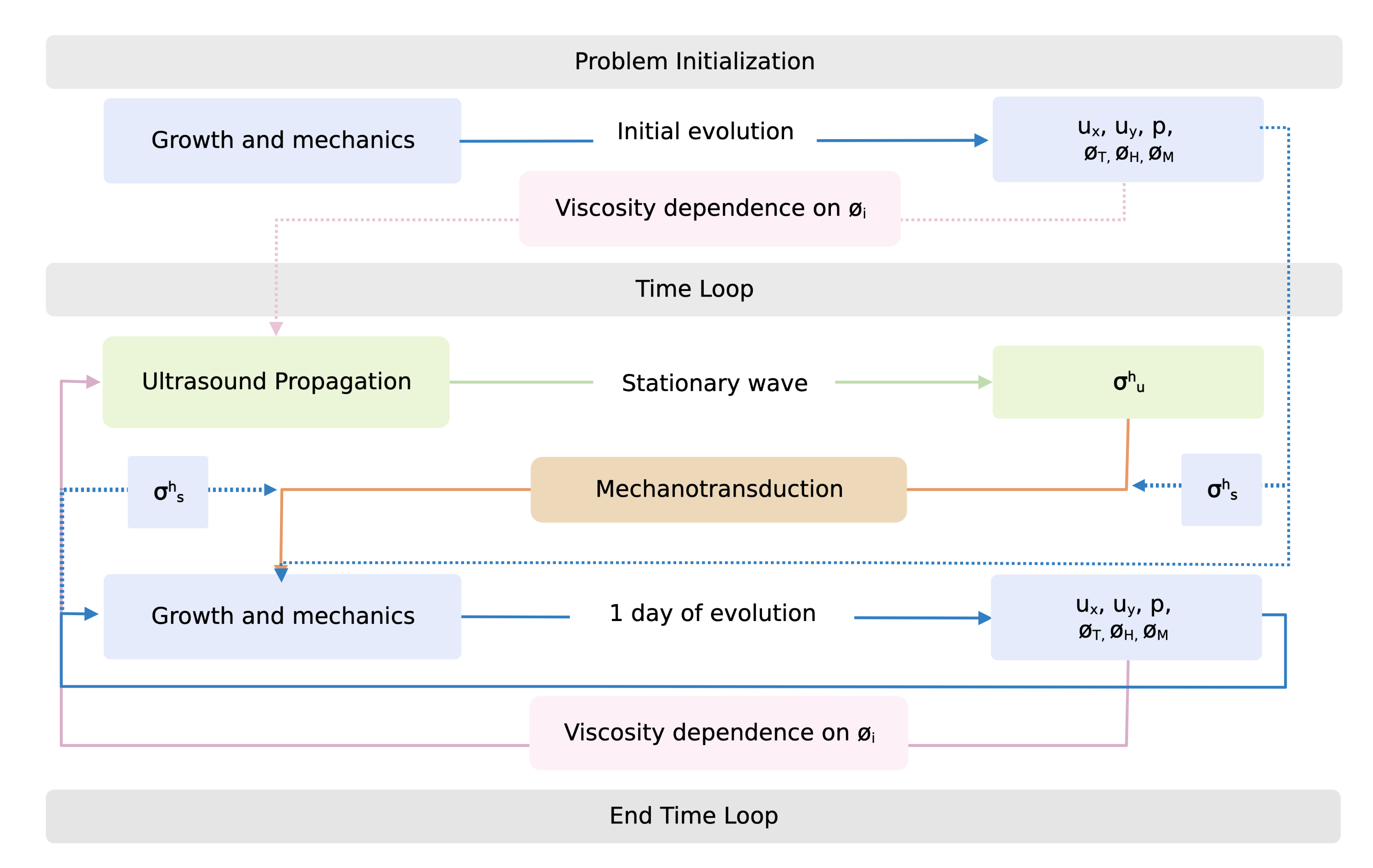}
\caption[Flowchart of the multiscale system]{\textbf{Flowchart of the multiscale system.} The system is initialized on a slow scale, where displacements, fluid pressure and solid phases are obtained. The solid phases are added to the fast-scale model of wave propagation to consider the viscosity of the tumor cell phase dependence, and it evolves until the stationary wave is achieved, where the ultrasonic hydrostatic stress is computed. Together with the slow hydrostatic stress, the ultrasonic stress is considered to compute the evolution of the system accounting for mechanotransduction. The results are again included at the slow and fast scale to complete the time loop until the final time of the simulations is achieved.}
\label{flowchart}
\end{figure*}

The self-developed computational model is solved in the Finite Element Analysis Program~\citep{taylor_feap_2014} -- FEAP --  and Matlab (MathWorks Inc., Natick, MA, USA), and visualized using Paraview~\citep{ParaView}.  We solve the multiscale system assuming a two-dimensional problem and plane strain, while the flowchart of the numerical simulations is reported in Figure~\ref{flowchart}.

Furthermore, we propose modeling ultrasound propagation independently of its duty cycle, regardless of whether it is continuous or pulsed. By applying the mechanical wave throughout the cell growth process, we ensure a comprehensive analysis. This decision is motivated by the fact that the duty cycle operates on a millisecond timescale and is based on persistence. Consequently, once a mechanotransduction stimulus is applied, the corresponding cellular response persists for several seconds. Thus, even during the silent period of the duty cycle, signaling remains activated, enabling us to consider the stress field throughout the ultrasonic period. In addition, we compute the ultrasound wave until it reaches a stationary state, where stress levels remain constant regardless of the duration of the sonication.  Through this methodology, we can effectively explore the implications of ultrasonic stress on the tumor microenvironment, unlocking its immense prospective as a promising avenue for cancer therapy.

 The parameter values used in slow-scale simulations are summarized in Table~\ref{table2} and the specific parameters for fast-scale ultrasound propagation include a range of frequencies between [1-20]MHz, acoustic pressures between [0.1-5]kPa, and viscosities between $\mathrm{[0-10]Pa \cdot s}$. These frequencies and intensities ranges are well below those established by the FDA, which considerably minimizes the possibility of collateral damage (spatial-peak temporal-average intensity ISPTA$<$100$\mathrm{mW \cdot cm^{-2}}$, and mechanical index MI$<$1.9). In addition, we simplify the degrees of freedom of the system and adapt parameters in our experiment prediction.

\begin{table}[!t]
\caption{\label{table2}Model parameters used in the simulations. \\}
\centering
\scriptsize
\resizebox{\textwidth}{!}{
\begin{tabular}{|c|c|c|c|c|c|}
\hline
\textbf{Description}       &    \textbf{Symbol}        & \textbf{Data} & \textbf{Units} & \textbf{Reference}  \\
\hline
Young modulus              & E &        $ 8$      & kPa             & \cite{roose2003solid,netti2000role,carotenuto2021lyapunov} \\ 
\hline
Undrained Poisson ratio       & $\nu_u$       & $0.49999$      & [-]  &\cite{roose2003solid,netti2000role}  \\ 
\hline
Drained Poisson ratio                      & $\nu$                     & $0.4$      & [-]  &\cite{roose2003solid,netti2000role,lorenzo2019}  \\ 
\hline
Mass density                 & $\rho$     &      $1000$     & $\mathrm{kg \cdot m^{-3}}$          & water \\ 
\hline
Hydraulic conductivity            & $k_h$                            & $3.1 \cdot 10^{-14}$      & $\mathrm{m^{2} \cdot Pa^{-1} s^{-1}}$        & \cite{roose2003solid,netti2000role,stylianopoulos2013coevolution,jain2007effect,wu2014effect}\\ 
\hline
Biot coefficient                       &   $\alpha$                           & $ 9.91\cdot 10^{-1}$      & [-]              &  \cite{carotenuto2021lyapunov,mario2021isogeometric}  \\ 
\hline
Biot modulus                          &    $M $                      & $1.79$     & $\mathrm{MPa}$            & \cite{carotenuto2021lyapunov,mario2021isogeometric}  \\ 
\hline
Vessel conductivity          & $k_v$                          & $2.70 \cdot 10^{-8}$     & $\mathrm{Pa^{-1}\cdot s^{-1}}$          &\cite{jain2007effect,stylianopoulos2013coevolution,wu2013effect}  \\ 
\hline
Vessel pressure              & $p_v$                           & $3.33\cdot 10^{3}$      & $\mathrm{Pa}$             & \cite{stylianopoulos2013coevolution} \\ 
\hline
Reflection coefficient       & $\omega$                      & $9.00\cdot 10^{-1}  $    & [-]              &\cite{jain2007effect,wu2013effect}  \\ 
\hline
Interstitial osmotic pressure     & $\pi_v -\pi_e$              & $1.33\cdot 10^{3}$     & $\mathrm{Pa}$             & \cite{jain2007effect,wu2013effect} \\ 
\hline
Lymphatic conductivity       & $k_{lo}$                           &$ 9.98\cdot 10^{-8} $     &  $\mathrm{Pa^{-1}\cdot s^{-1}}$ & \cite{stylianopoulos2013coevolution} \\ 
\hline
Lymphatic vessel pressure    & $p_l$                                & $1.33\cdot 10^{2}$      & $\mathrm{Pa}$&\cite{carotenuto2018growth}  \\ 
\hline
Exchange coefficient         & $\alpha_{TT}$                          & $1.30$     & [-]               & \cite{carotenuto2021lyapunov} \\ 
\hline
Exchange coefficient         &  $\alpha_{TH}$                             &$ 1.00$      & [-]               & \cite{carotenuto2021lyapunov}   \\ 
\hline
Exchange coefficient         &  $\alpha_{TM}$                                & $1.00$      & [-]               &  \cite{carotenuto2021lyapunov}  \\ 
\hline
Exchange coefficient         &  $\alpha_{HH}$                             &$ 3.00$      & [-]               & \cite{carotenuto2021lyapunov}   \\ 
\hline
Exchange coefficient         &  $\alpha_{HT}$                               & $2.00$     & [-]               &   \cite{carotenuto2021lyapunov} \\ 
\hline
Exchange coefficient         &  $\alpha_{HM}$                            & $1.00$      & [-]               &  \cite{carotenuto2021lyapunov}  \\ 
\hline
ECM production from $\phi_T$              & $ \beta_{T}$                             & $5.79\cdot 10^{-7}$      & $\mathrm{s^{-1}}$         &  \cite{carotenuto2021lyapunov}  \\ 
\hline
ECM production from $\phi_H$                   & $\beta_{H}$                              & $1.16\cdot 10^{-6}$      & $\mathrm{s^{-1}}$              &  \cite{carotenuto2021lyapunov}  \\
\hline
ECM degradation from $\phi_T$                 & $\delta_M \alpha_{MT} $                        & $2.89\cdot 10^{-6}$      & $\mathrm{s^{-1}}$              & \cite{carotenuto2021lyapunov}   \\ 
\hline
ECM degradation from $\phi_H$                 & $\delta_M \alpha_{MH} $                       & $2.89\cdot 10^{-6}$       & $\mathrm{s^{-1}}$             &   \cite{carotenuto2021lyapunov} \\ 
\hline
Initial condition $\phi^{T}$      & $\phi^{T0}$                                 & $1.50\cdot 10^{-1}$       & [-]              &  \cite{carotenuto2021lyapunov}  \\ 
\hline
Initial condition $\phi^{H}$      & $\phi^{H0}$                               & $1.50\cdot 10^{-1}$     & [-]               &   \cite{carotenuto2021lyapunov} \\
\hline
Initial condition $\phi^{M}$      & $\phi^{M0}$                           & $4\cdot 10^{-1}$       & [-]               &  \cite{carotenuto2021lyapunov}  \\ 
\hline
Proliferation rate $\phi^{T}$     & $T_T$                                 & $1.26\cdot 10^{-5}$       & $\mathrm{s^{-1}}$              &  \cite{carotenuto2021lyapunov}  \\ 
\hline
Proliferation rate $\phi^{H}$     & $T_H$                             & $1.26\cdot 10^{-5}$       & $\mathrm{s^{-1}}$              &  \cite{carotenuto2021lyapunov}  \\
\hline
Common lower rate            &$q$           & $0.05$      &     [-]           &  \cite{carotenuto2021lyapunov}  \\ 
\hline
Mechanotransduction smoothness  & $\chi_\sigma$                    & $ -0.05$       & $\mathrm{Pa^{-1}}$             &   \cite{carotenuto2021lyapunov} \\ Dynamic stress coefficient   & $\beta_s$                    & $0.2$       & $[-]$             &   fitted \\ 
\hline
Tumoral threshold stress             & $\sigma_{L}$               & $1.2\cdot 10^{3}$      & $\mathrm{Pa}$        & \cite{carotenuto2021lyapunov,helmlinger1997,roose2003solid,cheng2009} \\ 
\hline
\end{tabular}
}
\end{table}

We first reduce the degrees of freedom of the mathematical approach to fit the experimental data and reconstruct the mechanotransduction parameters. For simplicity, we have assumed the absence of the extracellular matrix and healthy phases, and we only consider the coexistence of proliferating tumor cells and fluid within the tumor spheroid, which means that $D_T=\alpha_{TH}=\alpha_{TM}=\Gamma_F=\Gamma_H=\Gamma_{M}=\beta_T=\beta_H=0$. We have chosen specific mechanical parameters from the experiment, including a frequency of $f=5$MHz and an acoustic pressure of $A=$1.5kPa, while tumor and culture medium viscosity $\eta_T$=2Pa$\cdot$s, and $\eta_c$=0.05Pa$ \cdot$s respectively, are assumed from the literature ranges~\citep{rus2020viscosity}.

 To estimate the total number of cells, despite the lack of experimental cell count localization, we integrate the tumor phase over space at a given time, represented as $\int_{A} \phi_T (x,y,t) dA$. We calibrate the simulation parameters using data from the control experiment, which takes into account the observed deceleration of cell proliferation on the first day, attributed to the rearrangement and development of spheroid clusters. Therefore, the absence of significant differences between the control and sonication groups on the first day may be attributed to cellular reorganization. From the first day to the third day, a consistent and higher proliferation rate is maintained. Once the growth parameters are calibrated, we further adjust the mechanotransduction parameters to align with the experimental results.

\subsection{Experimental setup}
The experimental procedure involved a total of three replicates for each condition, namely the control group and the sonicated group. For each replication, three independent experiments are conducted. 

\begin{figure*}[h]
\centering
\includegraphics[width=\textwidth]{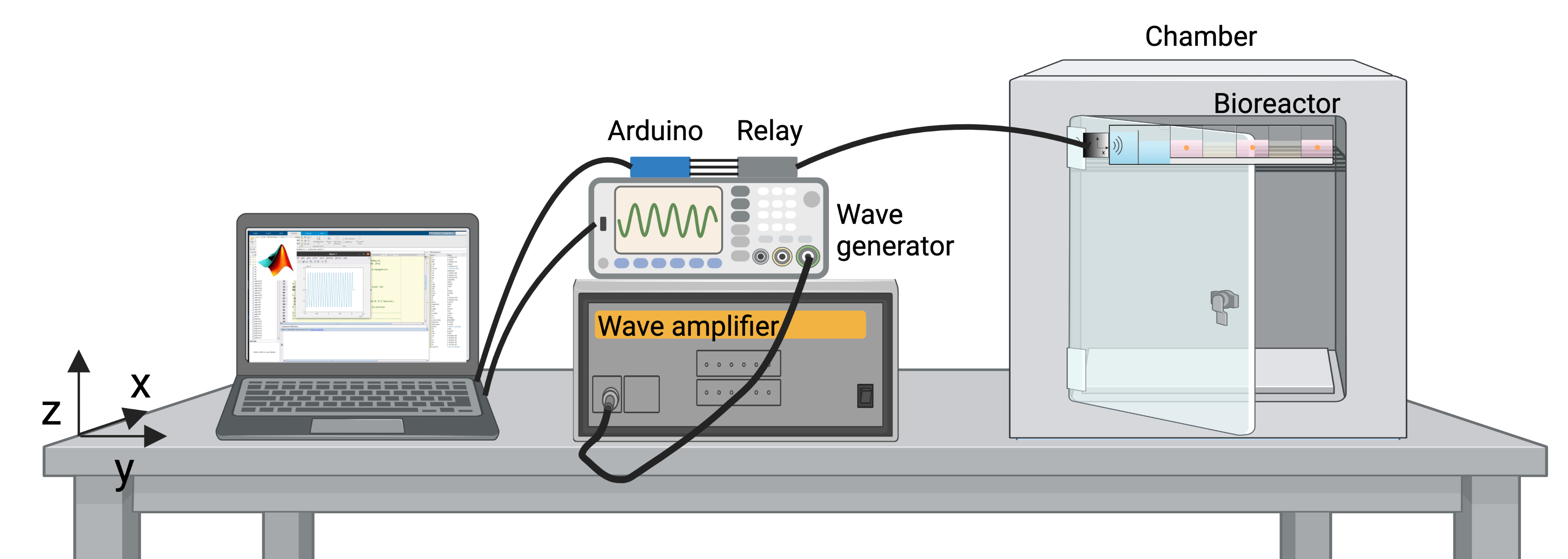}
\caption[Setup of the measurements]{\textbf{Setup of the measurements.} The Arduino is connected to the computer and the software is loaded, allowing the switch of mechanical signals. The Arduino also serves as a trigger to restart the signal and prevent any delay. The wave is generated using Matlab software and then loaded into the wave generator. Before connecting the wave generator to the amplifier, the signal is first verified using an oscilloscope to ensure that the frequencies and connections are correct.  Once the signals have been tested, the transducer is connected and placed on the support, and the coupling gel is extended on the transducers and bioreactor faces as a coupling material to avoid air bubbles. The relays are then connected. As the final step, cells are transferred to their designated chambers in the bioreactor and placed in the incubator until subsequent analysis.}
\label{setup_lab}
\end{figure*}

The bioreactor in which the CSCs are located is subjected to a 70$\%$ alcohol spray and subsequently placed in a chamber designed for ultraviolet (UV) sterilization. This process lasts 30 minutes and ensures disinfection of the bioreactor by effectively eliminating possible contaminants with the use of UV light. The water and attenuating material chambers are filled with caution to prevent overflow and avoid any form of contamination or interference. In addition, the chambers are covered with insulating tape. The complete setup is shown in Figure~\ref{setup_lab}. To conduct the examination of cell proliferation, the bioreactor is removed from the incubator and moved to a UV sterilization chamber with laminar flow to avoid contamination.

\begin{figure*}[h]
\centering
\includegraphics[width=1\textwidth]{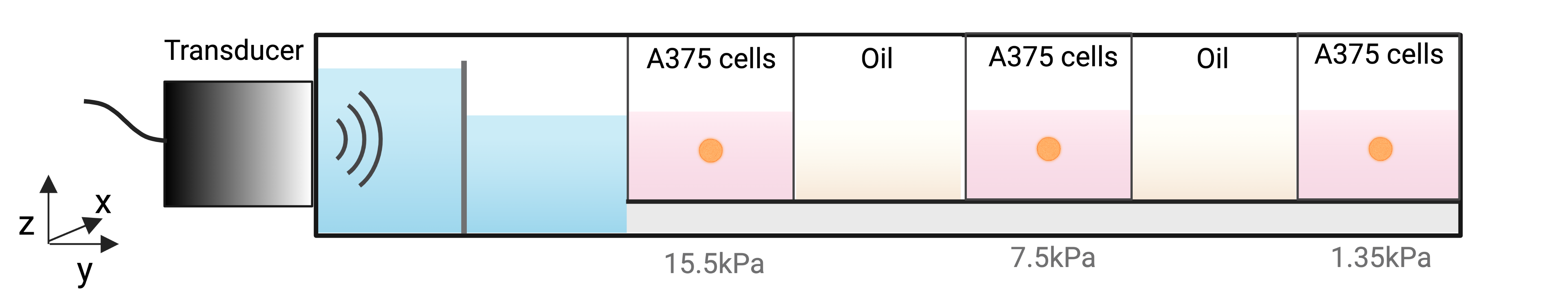}
\caption[Sonication scheme]{\textbf{Sonication scheme.} The transducer emits an ultrasonic wave through the first medium of water, which prevents the temperature from increasing. The wave then travels through the culture containing cells and attenuating media, causing the acoustic pressure to decrease as it encounters different materials and viscosities. As a result, the same bioreactor can be used for a given frequency and various acoustic pressures.}
\label{sonication_scheme}
\end{figure*}

The bioreactor used for cell culture consists of five sequentially arranged Petri dishes containing A-375 human melanoma cells embedded in a culture medium and an attenuating medium (oil), as depicted in Figure~\ref{sonication_scheme}. This experimental setup is designed to enable the generation of various wave amplitudes using a single transducer, as the emitted wave loses energy while propagating through different media. To prevent heating effects, a water-filled region is included at the beginning of the bioreactor. Acoustic pressure values are measured in each culture using a hydrophone probe, which is submerged in a replica of the bioreactor to capture acoustic pressure values without affecting tumor response. Through this method, we have determined that the first culture experiences 15.5 kPa, the second 7.5 kPa, and the third 1.35 kPa.

\subsection{Cell culture}
Melanoma cancer cell lines (A375) were acquired from the American Type Culture Collection (ATCC) and were cultured according to the procedures recommended by the ATCC. Cell lines were passaged for a period of fewer than 6 months and were regularly tested for mycoplasma contamination. Cells were maintained in advanced DMEM (Sigma-Aldrich) supplemented with 10$\%$ FBS (Gibco) and 5$ \%$ penicillin/streptomycin (Sigma-Aldrich).

To obtain tertiary spheres, melanoma cells were cultured in suspension using low-attachment plates containing DMEM-F12, 1$\%$ streptomy\-cin-penicillin, 1 mg/mL hydrocortisone (Sigma-Aldrich), 4 ng/mL heparin (Sigma-Aldrich), 1X ITS (Gibco), 1X B27 (Gibco), 10 ng/mL EGF (Sigma-Aldrich), 10 ng/mL FGF (Sigma-Aldrich), 10 ng/mL HGF and 10 ng/mL IL6 (Miltenyi Biotec) as previously described~\citep{lopez2022biofabrication}. Cells were cultured for 6 days and spheres are disaggregated every 72 hours until tertiary spheres were obtained. To achieve this, the spheres were collected by centrifugation at 1500 rpm for 5 minutes, incubated with trypsin-EDTA (Sigma-Aldrich) at 37 °C for 5 minutes, and then inactivated with FBS. The cells were then washed with PBS and reseeded under the same culture conditions. After that, melanospheres CSC phenotype  was confirmed as previously described~\citep{lopez2022biofabrication}.

\subsection{Cell proliferation assay}
Alamar Blue Assay (Biorad) was the measurement method. Cell growth was monitored on days 0, 1, and 3. To ensure reliable results, two parallel experiments were conducted. In the first experiment, a bioreactor loaded with melanoma CSCs was exposed to 24 hours of ultrasound and measurements were taken. In the second experiment, a bioreactor was used in which cells were treated for 72 hours without interruption. This approach was implemented to avoid any possible interference or damage during the manipulation of the spheroids.

The experimental protocol consisted of adding 10$\mu$l of Alamar Blue solution per 100$\mu$l of media to the cells and incubating them for 2 hours. Following the incubation period, the fluorescence intensity was measured using the Synergy HT instrument (BIO-TEK) at an excitation wavelength of 530nm and emissions of 590nm. For the data analysis, a non-parametric methodology was developed under the assumption of non-normality in the growth rate variables and the small size of the samples. The Kruskal-Wallis with Wilcox proves were performed for pairwise comparisons between group levels with corrections for multiple testing. RStudio software Version 1.4.1717 has been used to analyze the statistical differences. Although no apoptosis assays were performed in this study, we strongly encourage future research to include them to gain a more complete understanding of the underlying cellular processes.

\section{Results}
\subsection{LIUS hinders CSC growth \textit{in-vitro}}\label{sec2}
Melanoma CSCs are insonified at a frequency of 5MHz, enabling partial tumor penetration and enhanced mechanotransduction without cytoskeleton damage. Our experimental results yielded safe acoustic pressure values of 1.35kPa, 7.5kPa, and 15.5kPa, inducing mechanotransduction effectively without the tissue disruption associated with higher pressures.

In our experiments, we observed a significant decrease in the net proliferation of CSCs when subjected to 5MHz sonication, as compared to the control spheroids, over a three-day period (p=0.018, Wilcoxon-Mann-Whitney test). However, no significant differences were found among the different acoustic pressures. These results suggest that the sensitivity limit of cells may have been reached before reaching 1.5kPa, supporting the hypothesis that cells have a lower sensitivity limit for dynamic stress compared to static stress~\citep{carotenuto2021lyapunov,carotenuto2018growth,fraldi2018cells,helmlinger1997}.

The mathematical model used to explain experiments is described in Methods in the framework of the infinitesimal growth, poroelasticity theory~\citep{carotenuto2018growth,carotenuto2021lyapunov,roose2003solid,biot1941general} and a multiscale approach~\citep{kevorkian2012multiple,rus2014nature}.

\begin{figure*}[h]
\centering
\includegraphics[width=\linewidth]{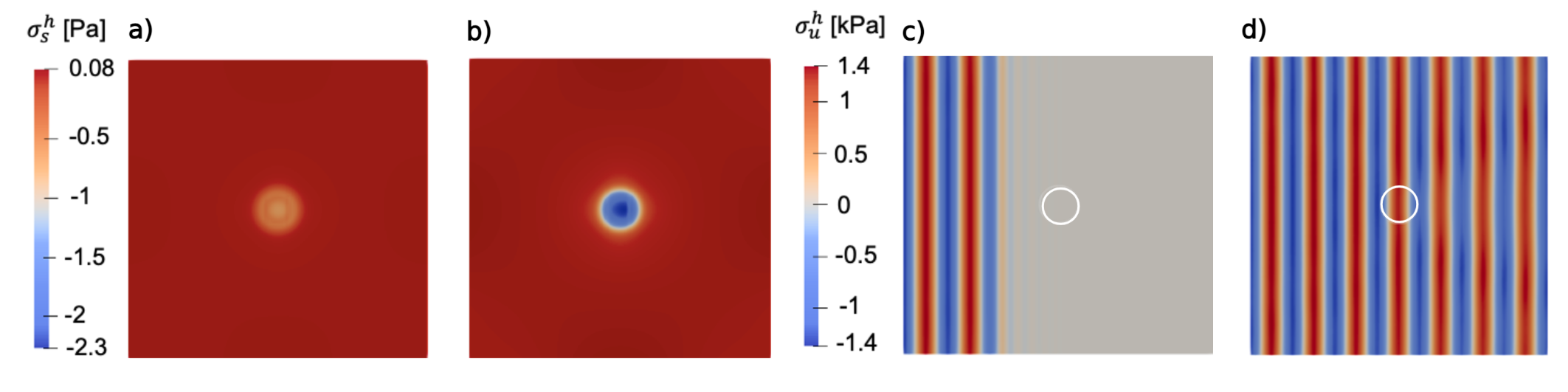}
\caption[Hydrostatic stresses during growth]{\textbf{Hydrostatic stresses during growth.} a) Slow hydrostatic stress of sonicated CSC at t = 1day  and b) slow hydrostatic stress of sonicated CSC at t = 3days. Slow hydrostatic compression increases at the center of the spheroid as it grows. c) Fast ultrasound stress before reaching tumor spheroid and d) fast ultrasound stress when stationary state is achieved, where a slight decrease in stress is perceived after reaching the tumor spheroid. The main parameters for sonication are: frequency $f=5$MHz, $A=1.5$kPa, tumor viscosity $\eta_T=2\mathrm{Pa\cdot s}$ and culture medium viscosity $\eta_c=0.05\mathrm{Pa\cdot s}$.}
\label{fig_stress}
\end{figure*}

From a computational standpoint, the hydrostatic stresses characterizing the stress state of the tumor manifest themselves at two distinct scales: the slow and ultrasonic stress,  as depicted in Figure~\ref{fig_stress}. With regard to slow-scale stress, compression is predominantly concentrated in the core of the tumor and consistent with previous research~\citep{northcott2018feeling,jain2014role,ramirez2017}. Additionally, this study also points out that the compression state increases over time and growth and operates at the order of Pascals, whereas ultrasound stress, which is defined by rarefaction and compression, is three orders of magnitude bigger. Such a substantial difference in stress scales implies that ultrasound stress is expected to exert a more prominent impact on mechanotransduction than the slow-scale stress. 

Furthermore, numerical simulations suggest that ultrasound diffraction through the tumor can result in slight shadow areas with lower displacements and stresses. For this particular case study, the slight diffraction presented is attributed to the difference in viscosity between the culture medium and the tumor spheroid at the applied frequency, resulting in a heterogeneous stress distribution within the bioreactor. 

\begin{figure*}[h]
\includegraphics[width=1\textwidth]{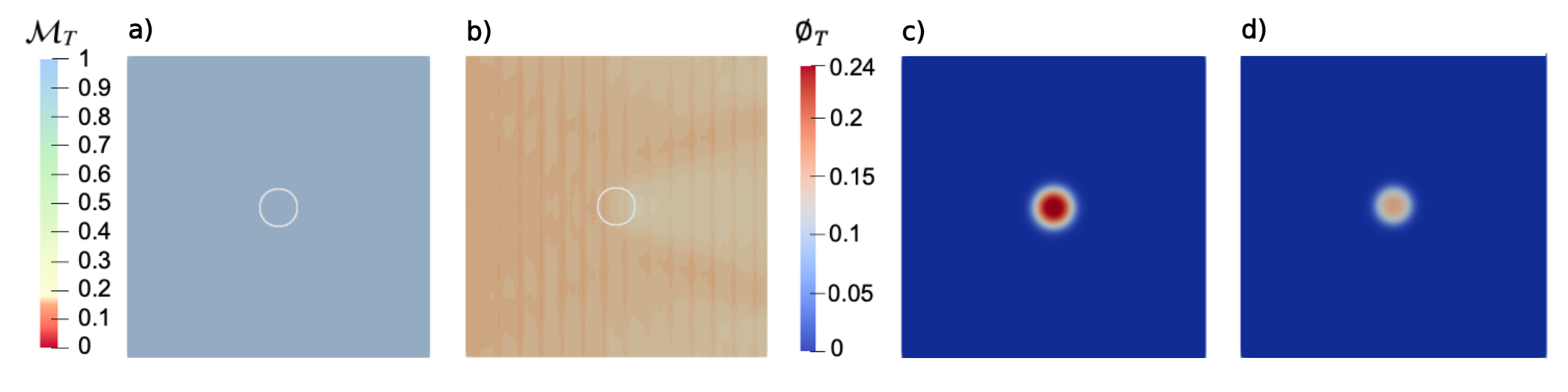}
\caption[Mechanotransduction and growth]{\textbf{Mechanotransduction and growth}. a) Control culture exhibits constant mechanotransduction values in space. b) Dynamic LIUS hydrostatic stress causes shadow areas that are translated into patterns in the mechanotransduction function. c) The control culture demonstrates a pronounced tumoral phase ($\phi_T$), representing a substantial fraction of tumor cells within the media composition. d) In comparison to the control cells, the sonicated cells exhibit attenuated proliferation, indicating a decrease in their growth rate. The main parameters used for these simulations are $\alpha_{TT}=2.9$, $\beta=0.2$, $t = 3$ days. The US parameters are $A=1.5$kPa, $f=5$MHz, $\eta_T=\mathrm{2Pa\cdot s}$ and $\eta_c=\mathrm{0.05Pa \cdot s}$.}
\label{fig2}
\end{figure*}

The propagation patterns of mechanical waves impact mechanotransduction processes. Computational analysis (Figure~\ref{fig2}) reveals that the mechanotransduction function $\mathcal{M}_T$ -- a mathematical relationship that describes how stress affects tumor growth -- remains spatially constant in the control culture. This is attributed to the total stress generated during the 3-day growth period remaining below the threshold stress of tumor cells for the duration of ultrasonic exposure, resulting in unaffected growth.
On the other hand, ultrasound diffraction introduces spatial heterogeneity in mechanotransduction. Our numerical observations suggest that the resulting stress shadow is inadequate to elicit discernible growth or migration patterns in cells via pressure gradients towards regions of lower stress. 

Therefore, the numerical results support a homogeneous decrease in tumor cell proliferation following sonication at a frequency of $f=5$ MHz relative to the control experiments, aligning with the proposed mechanotransduction mechanism.

\begin{figure*}[h]
\centering
\includegraphics[width=\textwidth]{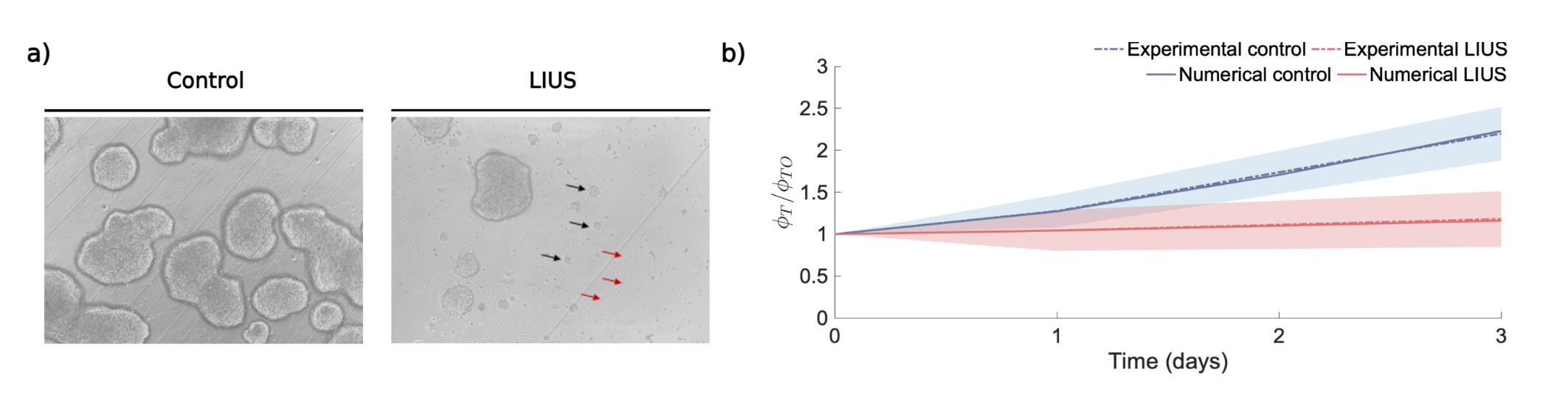}
\caption[Reduction in cell viability for sonicated spheroids]
{\textbf{Reduction in cell viability for sonicated spheroids at frequency $f=5$MHz, $A=1.5$kPa, tumor viscosity $\eta_T=2\mathrm{Pa\cdot s}$ and culture medium viscosity $\eta_c=0.05\mathrm{Pa\cdot s}$} a) Visualization of melanospheres with optical microscopy, where comparative analysis of untreated control and 5 MHz treated melanospheres after 72 hours is shown. Reduced-sized melanospheres are indicated by black arrows, while red arrows highlight individual cells. Images captured at 10x magnification. b) Computational model of LIUS mechanotherapy reproduces in vitro experiments. The dashed lines represent experiments, and the shaded bands their interval of confidence, while the solid lines denote numerical simulations. A trend change is observed between the control and sonicated cells, where cell proliferation decreases 48$\%$ on the third day after the application of LIUS. The proliferation parameters used are $T_T=0.58\cdot10^{-5}\mathrm{s^{-1}}$ on the first day and $T_T=0.77\cdot10^{-5} \mathrm{s^{-1}}$ from day one, while $\alpha_{TT}=2.9$. The mechanotransduction parameters are $q_T=0.05$, $b_T = 0.05$, $\beta_s=0.2$, and $\sigma_{L_T} = 1.2\mathrm{kPa}$.}
\label{fig_comparative}
\end{figure*}

Figure~\ref{fig_comparative}a clearly demonstrates that subjecting melanospheres to a 5 MHz frequency for 72 hours led to a significant decrease in both the number and size of the spheroids. Additionally, ultrasound had an adverse impact on cell viability, resulting in the identification of non-viable individual cells due to the induced disaggregation and toxic effects of LIUS.

The numerical simulations presented in this study closely replicate the initial experimental observations, as illustrated in Figure~\ref{fig_comparative}b. Despite the challenges of experimental cell count localization, our estimation method, based on integrating the tumor phase over space, provides insights into the overall cell population. By comparing the estimated cell count with experimental findings, we validate the accuracy of our simulation. Additionally, the calibration of simulation parameters using the control experiment data ensures the reliability of our results. The adjusted mechanotransduction parameters further enhance the agreement between our simulations and experimental outcomes.

\
In addition, we present a sensitivity analysis investigating the influence of frequencies, acoustic pressures, and viscosities on the behavior of tumor dynamics under various mechanical wave conditions.  

The key findings, summarized concisely in Table~\ref{table3}, provide valuable insights into the impact of different mechanical wave parameters on tumor proliferation rates. Our numerical simulations demonstrate that acoustic pressures greater than 1.5 kPa can produce substantial reductions in proliferation rates ranging from 46.6$\%$ to 48$\%$ for frequency values within the range of 1-5 MHz. Interestingly, these reductions occur when the viscosity of the medium is $\eta_c = \mathrm{0.05 Pa \cdot s}$, indicating that the perceived limitations of CSCs can be reached before surpassing this acoustic pressure threshold.
For higher frequency values, such as 20 MHz, an increase in acoustic pressure of 5 kPa is required to achieve a 52.5$\%$ decrease in proliferation rates. Furthermore, the influence of medium viscosity is observed, as higher viscosities result in increased wave attenuation, potentially limiting the reduction in proliferation rates and leading to values comparable to those of the control group (0-39.6$\%$ reduction compared to non-sonicated cells).

These findings and bioeffects align with previous research, as succinctly summarized in Table~\ref{tab_LIUS_exp}, where mechanotransduction studies typically involve frequencies ranging from 0.3 to 1.5 MHz. Notably, our study deviates with a higher frequency of 5 MHz, accompanied by substantially lower intensities around 5 kPa while achieving similar bioeffects. Furthermore, our parameters suggest that we are operating beyond the threshold of cytodisruption. Consequently, our results substantiate the hypothesis that comparable reductions in proliferation rates can be attained by applying lower levels of dynamic stress and higher frequencies, thereby reinforcing the effectiveness of these strategies for treating cancer stem cells. Additionally, we have reported a 50$\%$ reduction in proliferation rates compared to stress-free growth under static stress conditions, corroborating findings from prior studies~\citep{helmlinger1997,roose2003solid,montel2012,cheng2009,carotenuto2018growth,carotenuto2021lyapunov}.

\begin{table*}[!t]
\caption{\label{table3}Numerical results of LIUS therapy and proliferation reduction in comparison to the control group on the third day.\\}
\centering
\scriptsize
\begin{tabular}{|c|c|c|c|c|}
\hline
$f$ [MHz] & $A$ [kPa] & $\eta_c \mathrm{[Pa\cdot s]}$ & $\eta_T \mathrm{[Pa\cdot s]}$ & Proliferation decrease $[\%]$ \\
\hline
1   & $1.5$ & $0.05$ & $2$ & 48.4 \\ 
\hline
5   & $1.5$ & $0.05$ & $2$ & 48  \\
\hline
20  & $1.5$ & $0.05$ & $2$ & 27.3 \\
\hline
20  & $5$   & $0.05$ & $2$ & 52.5 \\
\hline
5    & $0.1$ & $0.05$ & $2$ & 0.0 \\
\hline
5    & $0.5$ & $0.05$ & $2$ & 38.1 \\
\hline
5   & $3$   & $0.05$ & $2$ & 50.7 \\
\hline
5   & $1.5$ & $0.05$ & $5$ & 47.5 \\
\hline
5   & $1.5$ & $0.05$ & $10$ & 46.6 \\
\hline
5   & $1.5$ & $2$     & $0.05$ & 39.6 \\
\hline
5   & $1.5$ & $5$     & $0.05$ & 14.8 \\
\hline
5   & $1.5$ & $10$    & $0.05$ & 0.21 \\
\hline
\end{tabular}
\end{table*}

\begin{table*}[!t]
\caption{\label{tab_LIUS_exp}Main setups in LIUS cancer therapy in vitro. The relation between intensity and acoustic pressure has been established. \\}
\centering
\resizebox{\textwidth}{!}{
\begin{tabular}{|p{2.5cm}|p{2cm}|p{2cm}|p{2cm}|p{4cm}|p{3cm}|p{1cm}|p{1.5cm}|}
\hline
\textbf{Cell line} & \textbf{Frequency $[$MHz$]$} & \textbf{Intensity $[$mW/cm$^2$$]$} & \textbf{Acoustic pressure $[$MPa$]$} & \textbf{Setup and Comments} & \textbf{Bioeffects} &\textbf{N}  & \textbf{Source}\\
\hline
CT-26, K562, U937, T cell  (in suspension)& 0.3-0.67& $ < 9.7 \cdot 10^{4}$ & $< 1.2$ &TUS= 2 min/day for 2 days. PD=2-40ms. DC= 10$\%$. Requires standing waves and reflection & Cytodisruption. Selective growth inhibition  &3-9& \cite{mittelstein2020selective} \\
\hline
T47D, MCF-12A (monolayers)  & 1.5 & 10,30,  \ \ \ 50,100  & 0.012, 0.021, 0.027, 0.039 & TUS=10min/day for 3 days. PD = 200$\mu$s. DC= 20$\%$. Decreasing proliferation with increasing intensity, PD, and DC. & Mechanotransduction. Selective growth inhibition &1 &\cite{katiyar2020inhibition} \\
\hline
HT29, Caco2  & 0.65-4.5  & 87.4 - $6.7 \cdot 10^{4}$ &0.036-1 &TUS= 10min/day for 1 day. PD= 30s. DC= 25$\%$ & Mechanotransduction - Cytodisruption. Growth inhibition  &2&~\cite{lucchetti2020low} \\
\hline
MDA-MB-231, Raw-264.7 &  1.5&  30   & 0.021 & TUS=20min/day for 10 days. PD = 200$\mu$. DC= 20$\%$  & Mechanotransduction. Reduction of osteoclastic differentiation & 3&~\cite{carina2018inhibitory} \\
\hline
A375, A549, Hela, Hacat &  0.67 &254  &0.061& TUS = 2min for 2 days, PD=30ms, DC= 10$\%$. Importance of stress field distribution. & Mechanotransduction - Cytodisruption. Selective growth inhibition  & $[-]$&~\cite{lin2022low} \\
\hline
MDA-MB-231, A375P, HT180 &  0.33 & 7.7  & 0.011 & TUS= 2h/day for 3 days. PD =  $[-]$. DC= 50$\%$ &Mechanotransduction: Piezo1 channel. Growth inhibition  & 2-3 &~\cite{tijore2020ultrasound}\\
\hline
MDA-MB-231, MCF10A (monolayers in matrigel) &  0.33 & 7.7 & 0.011 &TUS= 2h/day for 3 days. PD = $[-]$. DC= 50$\%$. Growth inhibition &Mechanotransduction: Piezo1 channel. Growth inhibition & 4 &~\cite{singh2021enhanced}\\
\hline
PANC-1 (monolayers)&  1 & $<100$ &  0.038 &TUS= 10-20-30min. DC= 100$\%$& Mechanotransduction. Migration inhibition& 4 &~\cite{gonzalez2023low}\\
\hline
\end{tabular}
}
\end{table*}

\subsection{LIUS causes selective patterns}\label{sec3}
To broaden the scope of our study, explore a wider range of phenomena, and gain a deeper understanding of the underlying mechanisms at play, we extend our research from the micro-scale to the macro-scale. Thus we propose to apply LIUS to a previously validated mathematical model, see~\citep{carotenuto2021lyapunov}. For these simulations, we use all the equations described in Methods and the mechanotransduction parameters fitted to our experimental data. We investigate the potential of 1MHz ultrasound as a selective modification tool for cancer cells with viscosities of $\mathrm{5Pa\cdot s}$.

\begin{figure*}[h]
\centering
\includegraphics[width=\textwidth]{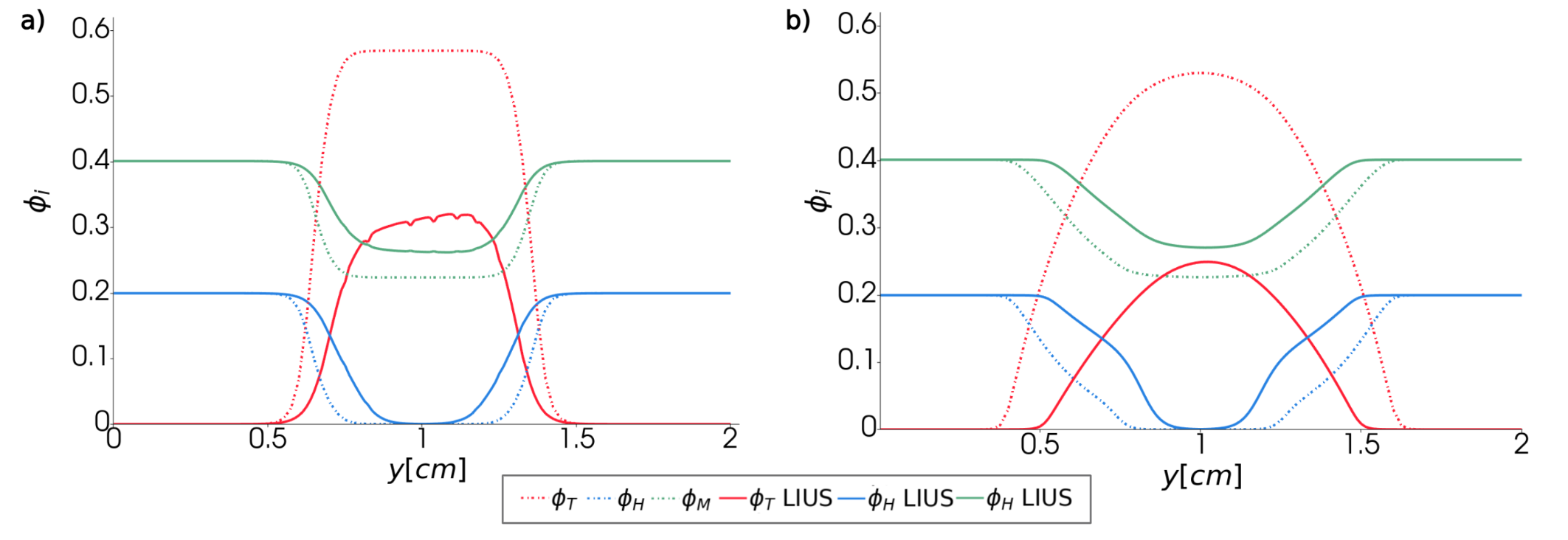}
\caption{\textbf{Patterns in growth and migration}. Dashed lines refer to tumor dynamics without sonication while solid lines refer to a sonicated tumor phase. a) When migration is not allowed, LIUS selectively reduces tumor cell phase proliferation, causing patterns in low-stress areas that translate to ECM phase while healthy phase remains unaltered. b) Tumor phase migration is reduced by LIUS sonication while it dampens patterns of low-stress areas, keeping the direction of tumor cell phase concentration. The used parameters are $f$=1MHz, $A$=1.5kPa, and $\eta_T=5$Pa$\cdot$s. Results at time $t=21$days. We refer to~\textit{Supplementary Material Video S2} for the animation.}
\label{patterns_selective}
\end{figure*}

Our simulations (Figure~\ref{patterns_selective}) reveal that ultrasound does not directly affect the proliferation and production of healthy cells. Stress experienced by healthy cells is considerably lower than that experienced by tumor cells, suggesting that ultrasound selectively attacks cancer cells, as reported in~\citep{mittelstein2020selective,katiyar2020inhibition,lin2022low}. We analyze the interplay of tumor and healthy cell growth, influenced by predator-prey dynamics. Our  results confirm a decrease in tumor cell growth that affects healthy cells without impacting their overall proliferation due to a higher proliferation threshold ($\sigma_L\geq10$kPa and $\beta_s=0.2$).

Furthermore, the existence of areas characterized by varying stress levels can trigger instabilities within tumors. This process disrupts their initial symmetry resulting in the concentration of cells in regions with low levels of stress and generating patterns as experimentally studied in~\citep{lin2022low}. The patterns are transmitted to the ECM phase, breaking into its homogeneous growth and possibly inducing its remodeling. 

The computational results also suggest that LIUS decreases tumor migration compared to nonsonicated cells, which is in agreement with the experimental results obtained in~\citep{gonzalez2023low}. Additionally, the results indicate that healthy cells and the extracellular matrix exhibit an adaptive growth response to cancer phase motility through a cross-diffusion phenomenon influenced by predator-prey dynamics, as it is shown in \textit{Supplementary Material Video S1}. Furthermore, migration dissipates and homogenizes differences in growth patterns  while keeping a predetermined direction of tumor cell concentration.

These observations provide a suitable explanation for the persistent proliferation of cells observed in certain experimental studies. Tumor cells may have the ability to proliferate in areas known as shadow zones, where the stress threshold is reached heterogeneously, resulting in a total cell count comparable to that of the control group. Considering motility, leader cancer cells located in low stress areas have the capacity to enhance migration, probably through the use of their cytonemes that sense and respond to the stress state of the surrounding environment~\citep{aguirre2022predictive,blanco2021modeling}, allowing them to retract from high stress areas. 

\begin{figure*}[h]
\includegraphics[width=1\textwidth]{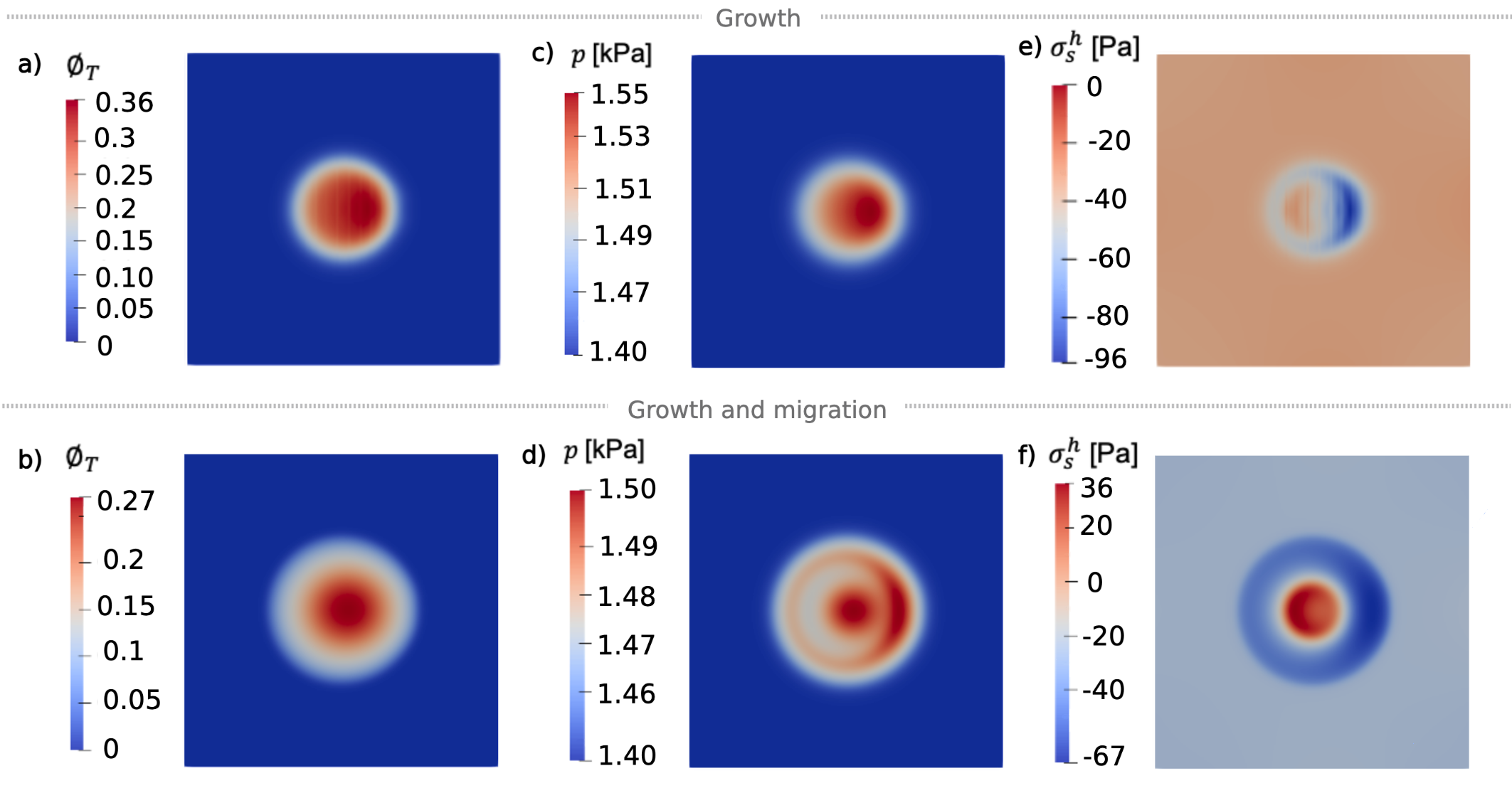}
   \caption[Tumor phase, interstial fluid pressure, and slowstress after LIUS sonication]{\textbf{Tumor phase, interstial fluid pressure, and slow stress after LIUS sonication}. Tumor phase, IFP and slow stress increase in the direction of wave propagation, while migration homogenizes the response of the tumor dynamics.  The main parameters used for these simulations are $\eta_T=10\mathrm{Pa\cdot s}$, $f$=1MHz, $A$=1.5kPa. Results are shown for t = 21 days. We refer to~\textit{Supplementary Material Video S3} for the animation.} 
\label{patterns_growth_excesive}
\end{figure*}

In fact, if tumor viscosity increases (see Figure~\ref{patterns_growth_excesive}), cells tend to concentrate in areas of lower stress, displacing the center of the tumor phase in the direction of wave propagation. Furthermore, the presence of cell mobility leads to a decrease in overall cell concentration but an increase in dispersion.

Interstitial fluid pressure increases during tumor growth, resulting in an uneven distribution of pressure. This pressure gradient could compress blood vessels, hindering the delivery of nutrients and drugs to tumor cells~\citep{stylianopoulos2013coevolution,stylianopoulos2018reengineering,blanco2023review}, while migration aids in stress dissipation and mitigates the substantial elevation of interstitial pressure observed when migration is disregarded \citep{blanco2021modeling}. 
Likewise, the growth of a sonicated tumor generates a slow and gradual increase in stress, although it is relatively minor compared to controlled growth conditions,
while compression occurs predominantly in the direction of tumor expansion, aligning with the preferred migration pathway.

\section{Discussion}

In this investigation, we unveil a theoretical framework to elucidate the potential mechanism through which LIUS selectively targets cancer stem cells. Merging computer simulations rooted in a mechanically coupled mathematical model of tumor spheroids with experimental validation, our results showcase that LIUS initiates a stress condition, strategically impeding and partially inhibiting tumor growth.

 Our numerical findings also demonstrate that ultrasound induces a significantly compressive hydrostatic stress state within the spheroid,  creating shadow areas due to viscosity differences between the medium and the tumor spheroid.

The proposed model considers the growth and migration of a poroelastic tumor, and the selection of strategic parameters based on the feedback with our experimental results. This allows us to address the possible challenges related to the prediction, controllability, and guidance of new experiments, while effectively avoiding effects associated with cavitation, cytodisruption and standing waves~\citep{mittelstein2020selective,heyden2016oncotripsy,heyden2017investigation,lin2022low,lucchetti2020low,prentice2005membrane}.  Regarding the complex interplay between LIUS and cancer dynamics, our study confirm that ultrasound does not directly affect healthy cell proliferation and production. Instead, ultrasound selectively targets cancer cells, minimizing adverse effects on surrounding healthy tissue. We show that reductions in proliferation rates can be achieved by applying lower levels of dynamic stress, reinforcing the potential of these strategies for cancer stem cell treatment. Likewise, our data allow us to conclude that LIUS decreases tumor migration, generating specific growth patterns. An increase in tumor viscosity or frequency leads to greater wave attenuation, resulting in wave diffraction that creates shadow areas with minor displacements and stresses, where the tumor cells are more prone to concentrating, growing, and migrating. These differences may help explain why, at times, the total cell count does not decrease compared to non-sonicated cells. 

Then, our interdisciplinary study provides a promising approach to exploring the effects of LIUS mechanotherapy on cancer stem cells, showing agreement with the previously known experimental results, and with the experiments developed for this paper.  Finally, additional validation using both longitudinal and shear waves~\citep{glatz2019low,hoelzlimpacting,blanco2023review} is necessary to confirm the selectivity of the attack. Once validated and refined, the model has the potential to bring about a significant transformation in the clinical approach to cancer treatment. Then, our results open up new prospects for further development and experimentation, and this LIUS strategy could offer a less aggressive effects on surrounding healthy tissue, more effective, and cost-effective treatment option for cancer stem cells. In the long run, integrating patient information, big data and artificial intelligence~\citep{lorenzo2023patient} holds the promise of tailoring LIUS-based treatments to individual patients, thus optimizing their effectiveness.

\section*{Acknowledgments}This study was funded by Ministry of Science, Innovation and Universities of Spain, project numbers PID2020-115372RB-I00 (B.B., M.H. and G.R.), PID2022-137228OB-I00 (J.S.), Consejería de Innovación, Ciencia y Empresa, Junta de Andalucía A-CTS-180-UGR20 (G.J, C.G and J.A.M.), B-FQM-580-UGR20 (B.B., J.S.), and by Consejería de Universidad, Investigación e Innovación from Junta de Andalucía, P21.\-00182 (B.B., M.H. and G.R.). This paper has been partially supported by the MINECO-FEDER (Spain) research grant number EQC2021-006920-P (J.A.M.) and PID2019-106947RA-C22 (B.B., J.M. and G.R.) and from the Chair 'Doctors Galera-Requena in cancer stem cell research' (CMC-CTS963). G.J. and C.G. acknowledge the posdoctoral fellowship from Plan Andaluz de Investigación, Desarrollo e Innovación (PAIDI 2020—FEDER funds—). Lastly, B.B. research was supported by the Ministry of Science, Innovation and Universities of Spain, FPU17\-01415. Under grant 101096884, Listen2Future is co-funded by the European Union. Views and opinions expressed are however those of the author(s) only and do not necessarily reflect those of the European Union or Key Digital Technologies Joint Undertaking. Neither the European Union nor the granting authority can be held responsible for them. The project is supported by the Key Digital Technologies Joint Undertaking and its members including top-up funding by Austria, Belgium, Czech Republic, Germany, Netherlands, Norway, and Spain. Figures of SI were created with Biorender.com.

\section*{Code availability}The codes that support the plots within this paper are described in the Methods, and they are available from the corresponding authors upon request.

\section*{Authors' contributions}The numerical model was self-coded in FEAP 8.6 by B. Blanco and R. Palma, while G. Rus contributed to the ultrasound propagation and mechanotransduction modeling in Matlab, and H. Gomez to weak formulations and code implementation. B. Blanco conducted the numerical simulations.
Concept and fundamentals of equations by B. Blanco, H. Gomez, R. Palma, G.Rus, and J. Soler.  Biological experiments were performed by G. Jiménez, C. Griñán-Lisón and J.A. Marchal while M. Hurtado, J. Melchor, G. Rus, and J. Soler designed the bioreactor and the experimental setup. The manuscript was reviewed for all authors.


\bibliographystyle{unsrtnat}
\bibliography{MechLIUS} 

\begin{thebibliography}{62}
\providecommand{\natexlab}[1]{#1}
\providecommand{\url}[1]{\texttt{#1}}
\expandafter\ifx\csname urlstyle\endcsname\relax
  \providecommand{\doi}[1]{doi: #1}\else
  \providecommand{\doi}{doi: \begingroup \urlstyle{rm}\Url}\fi

\bibitem[Huang and Kidoaki(2020)]{huang2020stiffness}
Daoxiang Huang and Satoru Kidoaki.
\newblock Stiffness-optimized drug-loaded matrix for selective capture and
  elimination of cancer cells.
\newblock \emph{J. Drug Deliv. Sci. Technol.}, 55:\penalty0 101414, 2020.
\newblock URL \url{https://doi.org/10.1016/j.jddst.2019.101414}.

\bibitem[Jain and Stylianopoulos(2010)]{jain2010delivering}
Rakesh~K Jain and Triantafyllos Stylianopoulos.
\newblock Delivering nanomedicine to solid tumors.
\newblock \emph{Nat. Rev. Clin. Oncol.}, 7\penalty0 (11):\penalty0 653, 2010.
\newblock URL \url{https://doi.org/10.1038/nrclinonc.2010.139}.

\bibitem[Polydorou et~al.(2017)Polydorou, Mpekris, Papageorgis, Voutouri, and
  Stylianopoulos]{polydorou2017pirfenidone}
Christiana Polydorou, Fotios Mpekris, Panagiotis Papageorgis, Chrysovalantis
  Voutouri, and Triantafyllos Stylianopoulos.
\newblock Pirfenidone normalizes the tumor microenvironment to improve
  chemotherapy.
\newblock \emph{Oncotarget}, 8\penalty0 (15):\penalty0 24506, 2017.
\newblock URL \url{https://doi.org/10.18632/oncotarget.15534}.

\bibitem[Panagi et~al.(2022)Panagi, Mpekris, Chen, Voutouri, Nakagawa, Martin,
  Hiroi, Hashimoto, Demetriou, Pierides, et~al.]{panagi2022polymeric}
Myrofora Panagi, Fotios Mpekris, Pengwen Chen, Chrysovalantis Voutouri,
  Yasuhiro Nakagawa, John~D Martin, Tetsuro Hiroi, Hiroko Hashimoto, Philippos
  Demetriou, Chryso Pierides, et~al.
\newblock Polymeric micelles effectively reprogram the tumor microenvironment
  to potentiate nano-immunotherapy in mouse breast cancer models.
\newblock \emph{Nat. Comm.}, 13\penalty0 (1):\penalty0 7165, 2022.

\bibitem[Abedi et~al.(2022)Abedi, Yao, Mittelstein, Bar-Zion, Swift,
  Lee-Gosselin, Barturen-Larrea, Buss, and Shapiro]{abedi2022ultrasound}
Mohamad~H Abedi, Michael~S Yao, David~R Mittelstein, Avinoam Bar-Zion,
  Margaret~B Swift, Audrey Lee-Gosselin, Pierina Barturen-Larrea, Marjorie~T
  Buss, and Mikhail~G Shapiro.
\newblock Ultrasound-controllable engineered bacteria for cancer immunotherapy.
\newblock \emph{Nat. Comm.}, 13\penalty0 (1):\penalty0 1585, 2022.

\bibitem[Mittelstein et~al.(2020)Mittelstein, Ye, Schibber, Roychoudhury,
  Martinez, Fekrazad, Ortiz, Lee, Shapiro, and
  Gharib]{mittelstein2020selective}
David~R Mittelstein, Jian Ye, Erika~F Schibber, Ankita Roychoudhury,
  Leyre~Troyas Martinez, M~Houman Fekrazad, Michael Ortiz, Peter~P Lee,
  Mikhail~G Shapiro, and Morteza Gharib.
\newblock Selective ablation of cancer cells with low intensity pulsed
  ultrasound.
\newblock \emph{Appl. Phys. Lett.}, 116\penalty0 (1):\penalty0 013701, 2020.
\newblock URL \url{https://doi.org/10.1063/1.5128627}.

\bibitem[Heyden and Ortiz(2016)]{heyden2016oncotripsy}
Stefanie Heyden and Michael Ortiz.
\newblock Oncotripsy: Targeting cancer cells selectively via resonant harmonic
  excitation.
\newblock \emph{J. Mech. Phys. Solids}, 92:\penalty0 164--175, 2016.
\newblock URL \url{https://doi.org/10.1016/j.jmps.2016.04.016}.

\bibitem[Heyden and Ortiz(2017)]{heyden2017investigation}
S~Heyden and M~Ortiz.
\newblock Investigation of the influence of viscoelasticity on oncotripsy.
\newblock \emph{Comput. Meth. Appl. Mech. Eng.}, 314:\penalty0 314--322, 2017.
\newblock URL \url{https://doi.org/10.1016/j.cma.2016.08.026}.

\bibitem[Lin et~al.(2022)Lin, Dong, Peng, Liu, Zhang, Lv, Yu, and
  Yao]{lin2022low}
Jianhao Lin, Shoulong Dong, Wencheng Peng, Hongmmei Liu, Penghao Zhang,
  Haoxiang Lv, Liang Yu, and Chenguo Yao.
\newblock Low-intensity pulsed ultrasound for killing tumor cells: The physical
  and biological mechanism.
\newblock pages 812--820, 2022.
\newblock URL \url{https://doi.org/10.1007/978-981-19-1528-4_83}.

\bibitem[Lucchetti et~al.(2020)Lucchetti, Perelli, Colella, Ricciardi-Tenore,
  Scoarughi, Barbato, Boninsegna, De~Maria, and Sgambato]{lucchetti2020low}
Donatella Lucchetti, Luigi Perelli, Filomena Colella, Claudio Ricciardi-Tenore,
  Gian~Luca Scoarughi, Gaetano Barbato, Alma Boninsegna, Ruggero De~Maria, and
  Alessandro Sgambato.
\newblock Low-intensity pulsed ultrasound affects growth, differentiation,
  migration, and epithelial-to-mesenchymal transition of colorectal cancer
  cells.
\newblock \emph{J. Cell. Physiol.}, 235\penalty0 (6):\penalty0 5363--5377,
  2020.
\newblock URL \url{https://doi.org/10.1002/jcp.29423}.

\bibitem[Prentice et~al.(2005)Prentice, Cuschieri, Dholakia, Prausnitz, and
  Campbell]{prentice2005membrane}
Paul Prentice, Alfred Cuschieri, Kishan Dholakia, Mark Prausnitz, and Paul
  Campbell.
\newblock Membrane disruption by optically controlled microbubble cavitation.
\newblock \emph{Nat. Phys.}, 1\penalty0 (2):\penalty0 107--110, 2005.

\bibitem[Katiyar et~al.(2020)Katiyar, Osborn, DasBanerjee, Zhang, Sarkar, and
  Sarker]{katiyar2020inhibition}
Amit Katiyar, Jenna Osborn, Malaya DasBanerjee, Lijie~Grace Zhang, Kausik
  Sarkar, and Krishna~Pada Sarker.
\newblock Inhibition of human breast cancer cell proliferation by low-intensity
  ultrasound stimulation.
\newblock \emph{J. Ultr. Med.}, 39\penalty0 (10):\penalty0 2043--2052, 2020.
\newblock URL \url{https://doi.org/10.1002/jum.15312}.

\bibitem[Carina et~al.(2018)Carina, Costa, Pagani, De~Luca, Raimondi, Bellavia,
  Setti, Fini, and Giavaresi]{carina2018inhibitory}
Valeria Carina, Viviana Costa, Stefania Pagani, Angela De~Luca, Lavinia
  Raimondi, Daniele Bellavia, Stefania Setti, Milena Fini, and Gianluca
  Giavaresi.
\newblock Inhibitory effects of low intensity pulsed ultrasound on
  osteoclastogenesis induced in vitro by breast cancer cells.
\newblock \emph{J. Exper. Clin. Cancer Res.}, 37\penalty0 (1):\penalty0 1--11,
  2018.
\newblock URL \url{https://doi.org/10.1186/s13046-018-0868-2}.

\bibitem[Gonz{\'a}lez et~al.(2023)Gonz{\'a}lez, Luzuriaga, Valdivieso, Frutos,
  L{\'o}pez, Hern{\'a}ndez, Rodr{\'\i}guez-Lorenzo, Yag{\"u}e, Santiago~Blanco,
  Pinto, et~al.]{gonzalez2023low}
Itziar Gonz{\'a}lez, Jon Luzuriaga, Alba Valdivieso, Jes{\'u}s Frutos, Jaime
  L{\'o}pez, Luis Hern{\'a}ndez, Luis~M Rodr{\'\i}guez-Lorenzo, Virginia
  Yag{\"u}e, Jos{\'e}~Luis Santiago~Blanco, Alberto Pinto, et~al.
\newblock Low-intensity continuous ultrasound to inhibit cancer cell migration.
\newblock \emph{Front. Cell Dev. Biol.}, 2023.
\newblock URL \url{https://doi.org/10.3389/fcell.2022.842965}.

\bibitem[Singh et~al.(2021)Singh, Tijore, Margadant, Simpson, Chitkara, Low,
  and Sheetz]{singh2021enhanced}
Aditi Singh, Ajay Tijore, Felix Margadant, Chloe Simpson, Deepak Chitkara,
  Boon~Chuan Low, and Michael Sheetz.
\newblock Enhanced tumor cell killing by ultrasound after microtubule
  depolymerization.
\newblock \emph{BioEng. Transl. Med.}, 6\penalty0 (3):\penalty0 e10233, 2021.
\newblock URL \url{https://doi.org/10.1002/btm2.10233}.

\bibitem[Tijore et~al.(2020)Tijore, Margadant, Yao, Hariharan, Chew, Powell,
  Bonney, and Sheetz]{tijore2020ultrasound}
Ajay Tijore, Felix Margadant, Mingxi Yao, Anushya Hariharan, Claire
  Alexandra~Zhen Chew, Simon Powell, Glenn~Kunnath Bonney, and Michael Sheetz.
\newblock Ultrasound-mediated mechanical forces selectively kill tumor cells.
\newblock \emph{BioRxiv}, 2020.

\bibitem[Na et~al.(2008)Na, Collin, Chowdhury, Tay, Ouyang, Wang, and
  Wang]{na2008rapid}
Sungsoo Na, Olivier Collin, Farhan Chowdhury, Bernard Tay, Mingxing Ouyang,
  Yingxiao Wang, and Ning Wang.
\newblock Rapid signal transduction in living cells is a unique feature of
  mechanotransduction.
\newblock \emph{Proc. Natl. Acad. Sci.}, 105\penalty0 (18):\penalty0
  6626--6631, 2008.
\newblock URL \url{https://doi.org/10.1073/pnas.0711704105}.

\bibitem[Geiger et~al.(2009)Geiger, Spatz, and
  Bershadsky]{geiger2009environmental}
Benjamin Geiger, Joachim~P Spatz, and Alexander~D Bershadsky.
\newblock Environmental sensing through focal adhesions.
\newblock \emph{Nat. Rev. Mol. Cell Biol.}, 10\penalty0 (1):\penalty0 21--33,
  2009.
\newblock URL \url{https://doi.org/10.1038/nrm2593}.

\bibitem[Vogel(2006)]{vogel2006mechanotransduction}
Viola Vogel.
\newblock Mechanotransduction involving multimodular proteins: converting force
  into biochemical signals.
\newblock \emph{Annu. Rev. Biophys. Biomol. Struct.}, 35:\penalty0 459--488,
  2006.
\newblock URL \url{https://doi.org/10.1146/Annu.Rev.biophys.35.040405.102013}.

\bibitem[Blanco et~al.(2023)Blanco, Gomez, Melchor, Palma, Soler, and
  Rus]{blanco2023review}
B~Blanco, H~Gomez, J~Melchor, R~Palma, J~Soler, and G~Rus.
\newblock Mechanotransduction in tumor dynamics modeling.
\newblock \emph{Phys. Life Rev.}, 44:\penalty0 279--301, 2023.
\newblock URL \url{https://doi.org/10.1016/j.plrev.2023.01.017}.

\bibitem[Broders-Bondon et~al.(2018)Broders-Bondon, Ho-Bouldoires,
  Fernandez-Sanchez, and Farge]{broders2018mechanotransduction}
Florence Broders-Bondon, Thanh Huong~Nguyen Ho-Bouldoires, Maria-Elena
  Fernandez-Sanchez, and Emmanuel Farge.
\newblock Mechanotransduction in tumor progression: the dark side of the force.
\newblock \emph{J. Cell Biol.}, 217\penalty0 (5):\penalty0 1571--1587, 2018.
\newblock URL \url{https://doi.org/10.1083/jcb.201701039}.

\bibitem[Agus and Michor(2012)]{agus2012sciences}
David~B Agus and Franziska Michor.
\newblock The sciences converge to fight cancer.
\newblock \emph{Nat. Phys.}, 8\penalty0 (11):\penalty0 773--774, 2012.

\bibitem[Colak and Medema(2014)]{colak2014cancer}
Selcuk Colak and Jan~P Medema.
\newblock Cancer stem cells--important players in tumor therapy resistance.
\newblock \emph{FEBS J.}, 281\penalty0 (21):\penalty0 4779--4791, 2014.

\bibitem[Olivares-Urbano et~al.(2020)Olivares-Urbano, Gri{\~n}{\'a}n-Lis{\'o}n,
  Marchal, and N{\'u}{\~n}ez]{olivares2020csc}
Mar{\'\i}a~Auxiliadora Olivares-Urbano, Carmen Gri{\~n}{\'a}n-Lis{\'o}n,
  Juan~Antonio Marchal, and Mar{\'\i}a~Isabel N{\'u}{\~n}ez.
\newblock {CSC} radioresistance: A therapeutic challenge to improve
  radiotherapy effectiveness in cancer.
\newblock \emph{Cells}, 9\penalty0 (7):\penalty0 1651, 2020.

\bibitem[Rus(2014)]{rus2014nature}
Guillermo Rus.
\newblock Nature of acoustic nonlinear radiation stress.
\newblock \emph{Appl. Phys. Lett.}, 105\penalty0 (12):\penalty0 121904, 2014.
\newblock URL \url{https://doi.org/10.1063/1.4894827}.

\bibitem[Biot(1941)]{biot1941general}
Maurice~A Biot.
\newblock General theory of three-dimensional consolidation.
\newblock \emph{J. Appl. Phys.}, 12\penalty0 (2):\penalty0 155--164, 1941.
\newblock URL \url{https://doi.org/10.1063/1.1712886}.

\bibitem[Carotenuto et~al.(2018)Carotenuto, Cutolo, Petrillo, Fusco, Arra,
  Sansone, Larobina, Cardoso, and Fraldi]{carotenuto2018growth}
AR~Carotenuto, A~Cutolo, A~Petrillo, R~Fusco, C~Arra, M~Sansone, D~Larobina,
  L~Cardoso, and M~Fraldi.
\newblock Growth and in vivo stresses traced through tumor mechanics enriched
  with predator-prey cells dynamics.
\newblock \emph{J. Mech. Behav. Biomed. Mater.}, 86:\penalty0 55--70, 2018.
\newblock URL \url{https://doi.org/10.1016/j.jmbbm.2018.06.011}.

\bibitem[Carotenuto et~al.(2021)Carotenuto, Cutolo, Palumbo, and
  Fraldi]{carotenuto2021lyapunov}
Angelo~Rosario Carotenuto, Arsenio Cutolo, Stefania Palumbo, and Massimiliano
  Fraldi.
\newblock Lyapunov stability of competitive cells dynamics in tumor
  mechanobiol.
\newblock \emph{Acta Mech. Sin.}, 37\penalty0 (2):\penalty0 244--263, 2021.
\newblock URL \url{https://doi.org/10.1007/s10409-021-01061-7}.

\bibitem[Kevorkian and Cole(2012)]{kevorkian2012multiple}
Jirair~K Kevorkian and Julian~D Cole.
\newblock \emph{Multiple scale and singular perturbation methods}, volume 114.
\newblock Springer Science \& Business Media, 2012.

\bibitem[Lorenzo et~al.(2019)Lorenzo, Hughes, Dominguez-Frojan, Reali, and
  Gomez]{lorenzo2019}
Guillermo Lorenzo, Thomas~JR Hughes, Pablo Dominguez-Frojan, Alessandro Reali,
  and Hector Gomez.
\newblock Computer simulations suggest that prostate enlargement due to benign
  prostatic hyperplasia mechanically impedes prostate cancer growth.
\newblock \emph{Proc. Natl. Acad. Sci.}, 116\penalty0 (4):\penalty0 1152--1161,
  2019.
\newblock URL \url{https://doi.org/10.1073/pnas.1815735116}.

\bibitem[Wu et~al.(2013)Wu, Frieboes, McDougall, Chaplain, Cristini, and
  Lowengrub]{wu2013effect}
Min Wu, Hermann~B Frieboes, Steven~R McDougall, Mark~AJ Chaplain, Vittorio
  Cristini, and John Lowengrub.
\newblock The effect of interstitial pressure on tumor growth: coupling with
  the blood and lymphatic vascular systems.
\newblock \emph{J. Theor. Biol.}, 320:\penalty0 131--151, 2013.

\bibitem[Fraldi and Carotenuto(2018)]{fraldi2018cells}
Massimiliano Fraldi and Angelo~R Carotenuto.
\newblock Cells competition in tumor growth poroelasticity.
\newblock \emph{J. Mech. Phys. Solids}, 112:\penalty0 345--367, 2018.
\newblock URL \url{https://doi.org/10.1016/j.jmps.2017.12.015}.

\bibitem[Stylianopoulos et~al.(2013)Stylianopoulos, Martin, Snuderl, Mpekris,
  Jain, and Jain]{stylianopoulos2013coevolution}
Triantafyllos Stylianopoulos, John~D Martin, Matija Snuderl, Fotios Mpekris,
  Saloni~R Jain, and Rakesh~K Jain.
\newblock Coevolution of solid stress and interstitial fluid pressure in tumors
  during progression: implications for vascular collapse.
\newblock \emph{Cancer Res.}, 73\penalty0 (13):\penalty0 3833--3841, 2013.
\newblock URL \url{https://doi.org/10.1158/0008-5472.can-12-4521f}.

\bibitem[Conte et~al.(2021)Conte, Casas-Tint{\`o}, and
  Soler]{conte2021modeling}
Martina Conte, Sergio Casas-Tint{\`o}, and Juan Soler.
\newblock Modeling invasion patterns in the glioblastoma battlefield.
\newblock \emph{PLoS Comp. Biol.}, 17\penalty0 (1):\penalty0 e1008632, 2021.
\newblock URL \url{https://doi.org/10.1371/journal.pcbi.1008632}.

\bibitem[Blanco et~al.(2021)Blanco, Campos, Melchor, and
  Soler]{blanco2021modeling}
Beatriz Blanco, Juan Campos, Juan Melchor, and Juan Soler.
\newblock Modeling interactions among migration, growth and pressure in tumor
  dynamics.
\newblock \emph{Math.}, 9\penalty0 (12):\penalty0 1376, 2021.
\newblock URL \url{https://doi.org/10.3390/math9121376}.

\bibitem[Calvo et~al.(2017)Calvo, Campos, Caselles, S{\'a}nchez, and
  Soler]{calvo2017qualitative}
Juan Calvo, Juan Campos, Vicent Caselles, {\'O}scar S{\'a}nchez, and Juan
  Soler.
\newblock Qualitative behaviour for flux-saturated mechanisms: travelling
  waves, waiting time and smoothing effects.
\newblock \emph{Journal of the European Mathematical Society}, 19\penalty0
  (2):\penalty0 441--472, 2017.
\newblock URL \url{http://dx.doi.org/DOI:10.4171/JEMS/670}.

\bibitem[Calvo et~al.(2016)Calvo, Campos, Caselles, S{\'a}nchez, and
  Soler]{calvo2016pattern}
Juan Calvo, Juan Campos, Vicent Caselles, O~S{\'a}nchez, and Juan Soler.
\newblock Pattern formation in a flux limited reaction--diffusion equation of
  porous media type.
\newblock \emph{Invent. Math.}, 206\penalty0 (1):\penalty0 57--108, 2016.
\newblock URL \url{https://doi.org/10.1007/s00222-016-0649-5}.

\bibitem[Hoffman et~al.(2011)Hoffman, Grashoff, and
  Schwartz]{hoffman2011dynamic}
Brenton~D Hoffman, Carsten Grashoff, and Martin~A Schwartz.
\newblock Dynamic molecular processes mediate cellular mechanotransduction.
\newblock \emph{Nature}, 475\penalty0 (7356):\penalty0 316--323, 2011.
\newblock URL \url{https://doi.org/10.1038/Nat.10316}.

\bibitem[Helmlinger et~al.(1997)Helmlinger, Netti, Lichtenbeld, Melder, and
  Jain]{helmlinger1997}
Gabriel Helmlinger, Paolo~A Netti, Hera~C Lichtenbeld, Robert~J Melder, and
  Rakesh~K Jain.
\newblock Solid stress inhibits the growth of multicellular tumor spheroids.
\newblock \emph{Nat. Biotechnol.}, 15\penalty0 (8):\penalty0 778, 1997.
\newblock URL \url{https://doi.org/10.1038/nbt0897-778}.

\bibitem[Roose et~al.(2003)Roose, Netti, Munn, Boucher, and
  Jain]{roose2003solid}
Tina Roose, Paolo~A Netti, Lance~L Munn, Yves Boucher, and Rakesh~K Jain.
\newblock Solid stress generated by spheroid growth estimated using a linear
  poroelasticity model.
\newblock \emph{Microvasc. Res.}, 66\penalty0 (3):\penalty0 204--212, 2003.
\newblock URL \url{https://doi.org/10.1016/s0026-2862(03)00057-8}.

\bibitem[Cheng et~al.(2009)Cheng, Tse, Jain, and Munn]{cheng2009}
Gang Cheng, Janet Tse, Rakesh~K Jain, and Lance~L Munn.
\newblock Micro-environmental mechanical stress controls tumor spheroid size
  and morphology by suppressing proliferation and inducing apoptosis in cancer
  cells.
\newblock \emph{PloS ONE}, 4\penalty0 (2):\penalty0 e4632, 2009.
\newblock URL \url{https://doi.org/10.1371/journal.pone.0004632}.

\bibitem[Dukhin and Goetz(2009)]{dukhin2009bulk}
Andrei~S Dukhin and Philip~J Goetz.
\newblock Bulk viscosity and compressibility measurement using acoustic
  spectroscopy.
\newblock \emph{J. Chem. Phys.}, 130\penalty0 (12):\penalty0 124519, 2009.
\newblock URL \url{https://doi.org/10.1063/1.3095471}.

\bibitem[Claes et~al.(2021)Claes, Chatwell, Baumh{\"o}gger, Hetk{\"a}mper,
  Zeipert, Vrabec, and Henning]{claes2021measurement}
Leander Claes, Ren{\'e}~Spencer Chatwell, Elmar Baumh{\"o}gger, Tim
  Hetk{\"a}mper, Henning Zeipert, Jadran Vrabec, and Bernd Henning.
\newblock Measurement procedure for acoustic absorption and bulk viscosity of
  liquids.
\newblock \emph{Measurement}, 184:\penalty0 109919, 2021.
\newblock URL \url{tps://doi.org/10.1016/j.measurement.2021.109919}.

\bibitem[d'Astous and Foster(1986)]{d1986frequency}
F.~T. d'Astous and F.~S. Foster.
\newblock Frequency dependence of ultrasound attenuation and backscatter in
  breast tissue.
\newblock \emph{Ultrasound Med. Biol.}, 12\penalty0 (10):\penalty0 795--808,
  1986.

\bibitem[Lysmer and Kuhlemeyer(1969)]{lysmer1969finite}
John Lysmer and Roger~L Kuhlemeyer.
\newblock Finite dynamic model for infinite media.
\newblock \emph{J. Eng. Mech.}, 95\penalty0 (4):\penalty0 859--877, 1969.

\bibitem[Taylor(2014)]{taylor_feap_2014}
R.~L. Taylor.
\newblock {FEAP} - finite element analysis program, 2014.
\newblock URL \url{http://www.ce.berkeley/feap}.

\bibitem[Ahrens et~al.(2005)Ahrens, Geveci, and Law]{ParaView}
James Ahrens, Berk Geveci, and Charles Law.
\newblock {ParaView}: An end-user tool for large data visualization.
\newblock In \emph{Visualization Handbook}. Elesvier, 2005.
\newblock {ISBN}~978-0123875822.

\bibitem[Netti et~al.(2000)Netti, Berk, Swartz, Grodzinsky, and
  Jain]{netti2000role}
Paolo~A Netti, David~A Berk, Melody~A Swartz, Alan~J Grodzinsky, and Rakesh~K
  Jain.
\newblock Role of extracellular matrix assembly in interstitial transport in
  solid tumors.
\newblock \emph{Cancer Res.}, 60\penalty0 (9):\penalty0 2497--2503, 2000.

\bibitem[Jain et~al.(2007)Jain, Tong, and Munn]{jain2007effect}
Rakesh~K Jain, Ricky~T Tong, and Lance~L Munn.
\newblock Effect of vascular normalization by antiangiogenic therapy on
  interstitial hypertension, peritumor edema, and lymphatic metastasis:
  insights from a mathematical model.
\newblock \emph{Cancer Res.}, 67\penalty0 (6):\penalty0 2729--2735, 2007.
\newblock URL \url{https://doi.org/10.1158/0008-5472.can-06-4102}.

\bibitem[Wu et~al.(2014)Wu, Frieboes, Chaplain, McDougall, Cristini, and
  Lowengrub]{wu2014effect}
Min Wu, Hermann~B Frieboes, Mark~AJ Chaplain, Steven~R McDougall, Vittorio
  Cristini, and John~S Lowengrub.
\newblock The effect of interstitial pressure on therapeutic agent transport:
  coupling with the tumor blood and lymphatic vascular systems.
\newblock \emph{J Theor Biol}, 355:\penalty0 194--207, 2014.
\newblock URL \url{https://doi.org/10.1016/j.jtbi.2014.04.012}.

\bibitem[de~Lucio et~al.(2021)de~Lucio, Bures, Ardekani, Vlachos, and
  Gomez]{mario2021isogeometric}
Mario de~Lucio, Miguel Bures, Arezoo~M Ardekani, Pavlos~P Vlachos, and Hector
  Gomez.
\newblock Isogeometric analysis of subcutaneous injection of monoclonal
  antibodies.
\newblock \emph{Comp. Meth. Appl. Mech. Eng.}, 373:\penalty0 113550, 2021.

\bibitem[Rus et~al.(2020)Rus, Faris, Torres, Callejas, and
  Melchor]{rus2020viscosity}
Guillermo Rus, Inas~H Faris, Jorge Torres, Antonio Callejas, and Juan Melchor.
\newblock Why are viscosity and nonlinearity bound to make an impact in
  clinical elastographic diagnosis?
\newblock \emph{Sensors}, 20\penalty0 (8):\penalty0 2379, 2020.
\newblock URL \url{https://doi.org/10.3390/s20082379}.

\bibitem[Lopez et~al.(2022)Lopez, Ruiz-Toranzo, Antich, Chocarro-Wrona,
  L{\'o}pez-Ru{\'\i}z, Jim{\'e}nez, and Marchal]{lopez2022biofabrication}
Julia Lopez, Marta Ruiz-Toranzo, Cristina Antich, Carlos Chocarro-Wrona, Elena
  L{\'o}pez-Ru{\'\i}z, Gema Jim{\'e}nez, and Juan~Antonio Marchal.
\newblock Biofabrication of a tri-layered 3d-bioprinted csc-based malignant
  melanoma model for personalized cancer treatment.
\newblock \emph{Biofabrication}, 2022.

\bibitem[Northcott et~al.(2018)Northcott, Dean, Mouw, and
  Weaver]{northcott2018feeling}
Josette~M Northcott, Ivory~S Dean, Janna~K Mouw, and Valerie~M Weaver.
\newblock Feeling stress: The mechanics of cancer progression and aggression.
\newblock \emph{Front. Cell Dev. Biol.}, 6:\penalty0 17, 2018.
\newblock URL \url{https://doi.org/10.3389/fcell.2018.00017}.

\bibitem[Jain et~al.(2014)Jain, Martin, and Stylianopoulos]{jain2014role}
Rakesh~K Jain, John~D Martin, and Triantafyllos Stylianopoulos.
\newblock The role of mechanical forces in tumor growth and therapy.
\newblock \emph{Annu. Rev. Biomed. Eng.}, 16:\penalty0 321--346, 2014.
\newblock URL \url{https://doi.org/10.1146/Annu.rev-bioeng-071813-105259}.

\bibitem[Ram{\'\i}rez-Torres et~al.(2017)Ram{\'\i}rez-Torres,
  Rodr{\'\i}guez-Ramos, Merodio, Penta, Bravo-Castillero, Guinovart-D{\'\i}az,
  Sabina, Garc{\'\i}a-Reimbert, Sevostianov, and Conci]{ramirez2017}
Ariel Ram{\'\i}rez-Torres, Reinaldo Rodr{\'\i}guez-Ramos, Jos{\'e} Merodio,
  Raimondo Penta, Juli{\'a}n Bravo-Castillero, Ra{\'u}l Guinovart-D{\'\i}az,
  Federico~J Sabina, Catherine Garc{\'\i}a-Reimbert, Igor Sevostianov, and Aura
  Conci.
\newblock The influence of anisotropic growth and geometry on the stress of
  solid tumors.
\newblock \emph{Int. J. Eng. Sci.}, 119:\penalty0 40--49, 2017.
\newblock URL \url{https://doi.org/10.1016/j.ijengsci.2017.06.011}.

\bibitem[Montel et~al.(2012)Montel, Delarue, Elgeti, Vignjevic, Cappello, and
  Prost]{montel2012}
Fabien Montel, Morgan Delarue, Jens Elgeti, Danijela Vignjevic, Giovanni
  Cappello, and Jacques Prost.
\newblock Isotropic stress reduces cell proliferation in tumor spheroids.
\newblock \emph{New J. Phys.}, 14\penalty0 (5):\penalty0 055008, 2012.
\newblock URL \url{https://doi.org/10.1088/1367-2630/14/5/055008}.

\bibitem[Aguirre-Tamaral et~al.(2022)Aguirre-Tamaral, Camb{\'o}n, Poyato,
  Soler, and blanco]{aguirre2022predictive}
Adri{\'a}n Aguirre-Tamaral, Manuel Camb{\'o}n, David Poyato, Juan Soler, and
  Isabel blanco.
\newblock {Predictive model for cytoneme guidance in Hedgehog signaling based
  on Ihog-Glypicans interaction}.
\newblock \emph{Nat. Comm.}, 13\penalty0 (1):\penalty0 1--14, 2022.
\newblock URL \url{https://doi.org/10.1038/s41467-022-33262-4}.

\bibitem[Stylianopoulos et~al.(2018)Stylianopoulos, Munn, and
  Jain]{stylianopoulos2018reengineering}
Triantafyllos Stylianopoulos, Lance~L Munn, and Rakesh~K Jain.
\newblock Reengineering the physical microenvironment of tumors to improve drug
  delivery and efficacy: from mathematical modeling to bench to bedside.
\newblock \emph{Trends Cancer}, 4\penalty0 (4):\penalty0 292--319, 2018.
\newblock URL \url{https://doi.org/10.1016/j.trecan.2018.02.005}.

\bibitem[Glatz(2019)]{glatz2019low}
Marlies Glatz.
\newblock Low frequency shear waves as a potential mechanotherapy approach in
  cancer.
\newblock 2019.

\bibitem[Hoelzl et~al.(20167)Hoelzl, Festy, Fruhwirth, and
  Sinkus]{hoelzlimpacting}
Marlies~Christina Hoelzl, Frederic Festy, Gilbert Fruhwirth, and Ralph Sinkus.
\newblock Impacting cancer cells via mechanical waves: can we change cellular
  behaviour?
\newblock In \emph{Proc. Int. Soc. Magn. Reson. Med.}, page 4363, 20167.

\bibitem[Lorenzo et~al.(2023)Lorenzo, Ahmed, Hormuth~II, Vaughn,
  Kalpathy-Cramer, Solorio, Yankeelov, and Gomez]{lorenzo2023patient}
Guillermo Lorenzo, Syed~Rakin Ahmed, David~A Hormuth~II, Brenna Vaughn,
  Jayashree Kalpathy-Cramer, Luis Solorio, Thomas~E Yankeelov, and Hector
  Gomez.
\newblock Patient-specific, mechanistic models of tumor growth incorporating
  artificial intelligence and big data.
\newblock \emph{arXiv preprint arXiv:2308.14925}, 2023.

\end{thebibliography}

\end{document}